\documentclass[10pt]{amsart}

\newcommand{\debug}{false}

\usepackage{ifthen}
\usepackage{enumerate}
\usepackage[dvips]{graphicx}
\usepackage{psfrag}
\usepackage{amsmath}
\usepackage{amssymb}

\newcommand{\pichere}[2]{
\ifthenelse{\equal{\debug}{true}}
{\begin{center}
{\fbox{\bf Figure \ref{fig:#2} Here}}
\end{center}}
{\includegraphics[width=#1\textwidth]{#2}}
}

\newcommand{\lab}[3]{\psfrag{#1}[#3]{$\scriptstyle{#2}$}}
\newcommand{\Lab}[3]{\psfrag{#1}[#3]{$\scriptscriptstyle{#2}$}}

\newtheorem{thm}{Theorem}
\newtheorem{lem}[thm]{Lemma}
\newtheorem{conj}[thm]{Conjecture}
\newtheorem{cor}[thm]{Corollary}

\theoremstyle{definition}
\newtheorem{defn}[thm]{Definition}
\newtheorem{defns}[thm]{Definitions}
\newtheorem{exmp}[thm]{Example}

\newenvironment{sketchproof}{\par\noindent{\sl Sketch Proof.}}{\qed}
\newenvironment{rem}{\par\noindent{\sl Remark.}}{\smallskip}
\newenvironment{rems}{\par\noindent{\sl Remarks.}}{\smallskip}

\newcommand{\cP}{{\mathcal{P}}}
\newcommand{\cD}{{\mathcal{D}}}
\newcommand{\cL}{{\mathcal{L}}}
\newcommand{\cO}{{\mathcal{O}}}
\newcommand{\cF}{{\mathcal{F}}}
\newcommand{\cA}{{\mathcal{A}}}

\newcommand{\LFP}{{\mathrm{LFP}}}
\newcommand{\RFP}{{\mathrm{RFP}}}
\newcommand{\LFS}{{\mathrm{LFS}}}

\newcommand{\bt}{{\mathrm{bt}}}
\newcommand{\BT}{{\mathrm{BT}}}

\newcommand{\ovun}[3]{\overset{#1}{\underset{#2}{#3}}}

\newcommand{\ov}[2]{\overset{#1}{#2}}
\newcommand{\un}[2]{\underset{#1}{#2}}

\newcommand{\card}{{\it \#}}
\newcommand{\I}{{^{-1}}}
\newcommand{\TTp}{TT$\,'$}
\newcommand{\TT}{TT\,}
\newcommand{\LD}{LD\,}
\newcommand{\be}{{\overline{e}}}
\newcommand{\opti}
{\raisebox{-.04cm}{$\stackrel{{\scriptscriptstyle 0}}{{\scriptscriptstyle
1}}$}}
\newcommand{\Rset}{{\mathbb R}}
\newcommand{\Zset}{{\mathbb Z}}
\newcommand{\Qset}{{\mathbb Q}}
\newcommand{\Nset}{{\mathbb N}}

\begin{document}
\bibliographystyle{amsalpha}

\title{Braid forcing and star-shaped train tracks}
\author{Andr\'e de Carvalho} 
\address{Institute for Mathematical Sciences, State University of New
York at Stony Brook, NY 11794-3660, USA}
\email{andre@math.sunysb.edu}
\author{Toby Hall}
\address{Department of Mathematical Sciences, University of
Liverpool, Liverpool L69 7ZL, UK}
\email{T.Hall@liv.ac.uk}
\begin{abstract}
Global results are proved about the way in which Boyland's forcing
partial order organizes a set of braid types: those of periodic orbits
of Smale's horseshoe map for which the associated train track is a
{\em star}. This is a special case of a conjecture introduced
in~\cite{dCH2}, which claims that forcing organizes all horseshoe
braid types into linearly ordered families which are, in turn,
parameterized by homoclinic orbits to the fixed point of code $0$.
\end{abstract}

\maketitle

\section{Introduction and background}

Given a discrete dynamical system and information about one of its
periodic orbits, can one derive further information about the system:
for example, the existence of other periodic orbits, or that it has
positive topological entropy?  The problem of periodic orbit forcing
in particular --- that is, to determine whether the presence of a
certain periodic orbit implies the existence of other periodic orbits
--- has interested dynamicists for many years. Sharkovskii's theorem
for self-maps of the real line is unrivalled in elegance: it defines a
total order $\preceq$ on the postive integers with the property that
if $m\preceq n$ then any continuous self-map of $\Rset$ which has a
periodic orbit of period $n$ must also have a periodic orbit of period
$m$. Moreover, the theorem is sharp: for any initial segment of the
order, there exists a continuous self-map of $\Rset$ having periodic
orbits of the corresponding periods and no others. In this context,
periodic orbits are specified by their period alone.  One could also,
for example, use the order on $\Rset$ to specify a periodic orbit by
its permutation, and then consider forcing among permutations. Many
authors have studied this problem, which is now quite well understood.

The corresponding problem in dimension~2 (i.e. for self-homeomorphisms
of the disk $D^2$) is much harder. In this case, the period alone is
an inadequate specification of a periodic orbit: given any set of
positive integers which includes~$1$, a self-homeomorphism of $D^2$
(or indeed of any other surface) can be constructed which has periodic
orbits of the given periods and no others. Boyland and others observed
that, by analogy with the permutation in dimension~1, adding
topological information about the way in which the points of a
periodic orbit of a disk homeomorphism `braid around' one another
produces a non-trivial theory. More precisely, Boyland defined the
{\em braid type} $\bt(P,f)$ of a periodic orbit $P$ of a disk
homeomorphism $f$ to be the isotopy class of $f$ relative to $P$, up
to topological change of coordinates. A braid type $\gamma$ is then
said to {\em force} another braid type $\beta$ if every disk
homeomorphism having a periodic orbit of type $\gamma$ also has one of
type $\beta$. Boyland showed that forcing is a partial order on the
set of braid types. This leads to the question of how forcing
organizes this set. Of course, if the question is to be answered, one
must be able to decide how two given braid types are related by the
partial order. This can be done using Bestvina and Handel's
algorithmic proof of Thurston's classification theorem for surface
homeomorphisms up to isotopy. However, not only is this process very
time consuming in practice, but also the information it gives is only
local, in the sense that the ability to compare two given braid types
gives no information about the global structure of the partially
ordered set of braid types.

This paper gives global information about the restriction of the
forcing order to a subset of braid types, namely those of periodic
orbits of Smale's horseshoe map for which the associated train track
is a {\em star} (that is, a tree with exactly one vertex of valence
larger than~1).  In~\cite{dCH2}, a conjecture is stated which
describes how forcing organizes all braid types of the horseshoe. It
claims that the symbolic code of each periodic orbit can be parsed
into two segments, the {\em prefix} and the {\em decoration}. All
orbits with the same decoration have the same {\em topological train
track type}, and form a family which is totally ordered by the forcing
relation; the position of an orbit within this totally ordered set is
determined by its prefix. Within families, this trivializes the
problem of comparing braid types: simply compare their symbolic codes
using the unimodal order. The conjecture also describes the forcing
relation between families in terms of the forcing between homoclinic
orbits associated to each family.

In this paper the conjecture is proved for the (infinite) family of
decorations for which the corresponding train track type is a
star. There is one such decoration for each rational
$m/n\in(0,1/2]$: the train track associated to $m/n$ is a star with
$n$ edges, and the train track map rotates the central vertex of the
star by $m/n$.

The approach taken is to start with train track maps of the
appropriate topological type, and to identify those of their periodic
orbits for which the star is itself a train track (after it has been
truncated outside the span of the orbit). The combinatorics of such
orbits are intricately related to the position of the rational $m/n$
within the Farey graph.

Only after these {\em train track orbits} have been determined is
their relationship with periodic orbits of Smale's horseshoe
investigated. Conceptually, the most important relationship between
the horseshoe and the star maps (and indeed the relationship which
motivates the conjecture under discussion) is that the star maps (or,
more accurately, the corresponding thick tree maps) can be obtained
from the horseshoe by {\em pruning}: that is, by performing an isotopy
which destroys all of the dynamics within a given open subset of the
disk, while leaving the dynamics unchanged elsewhere. Since this
construction is not used elsewhere in the paper, however, it is only
described on an intuitive level, and a less general approach to
showing that star periodic orbits have horseshoe braid types is
adopted. 

Section~\ref{sec:hsintro} describes the parsing of the code of a
horseshoe periodic orbit into prefix and decoration, and provides a
statement of the conjecture of~\cite{dCH2} and the special case of it
which will be proved here. Section~\ref{sec:farey} contains a brief
summary of the properties of the Farey graph which will be used in the
remainder of the paper. The main results of the paper can be found in
Section~\ref{sec:perorb}, in which star maps and the concept of a
train track orbit are defined, and the set of such orbits is
determined. Section~\ref{sec:prune} details the connection between the
star maps and the horseshoe, and Section~\ref{sec:app} is devoted to
the proof of a technical lemma.

The paper is often technical: it is the nature of the subject. Its
structure, however, is simple. Lemma~\ref{lem:identify} identifies the
periodic orbits of the horseshoe which have the property that their
period is equal to the sum of the denominators of the endpoints of
their rotation interval. The symbolic codes of these orbits have the
property that the prefix is determined by one endpoint of the rotation
interval, and the decoration by the other. Having described the
combinatorics of train track orbits of star maps
(Section~\ref{sec:perorb}), it is possible to compute their rotation
intervals, which have the property considered in
Lemma~\ref{lem:identify}. Since star orbits have horseshoe braid type
(Section~\ref{sec:prune}), it follows that the symbolic codes of the
train track orbits are as given by the lemma. The total order within
families follows directly from the nature of the maps being considered.

\subsection{Prefix and decoration for horseshoe periodic orbits}
\label{sec:hsintro}
\begin{figure}[htbp]
\lab{a}{0}{l}
\lab{b}{1}{l}
\lab{c}{S}{t}
\begin{center}
\pichere{0.3}{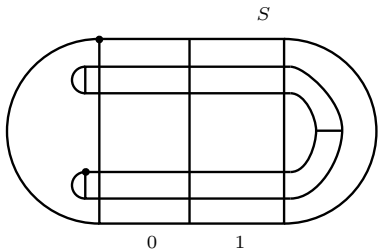}
\end{center}
\caption{Symbolic dynamics for the horseshoe}
\label{fig:hs}
\end{figure}

In this paper the standard model $F\colon D^2\to D^2$ of Smale's
horseshoe map~\cite{Smale} depicted in Fig.~\ref{fig:hs} will be
used; symbolic dynamics in the set $\Sigma_2=\{0,1\}^\Zset$ is applied
in the usual way to describe points $x\in D^2$ whose (past and future)
orbits lie entirely in the square $S$. The definitions and results
summarized in this section can be found in~\cite{Ha,dCH2}.

A periodic orbit $P$ of $F$ of (least) period $n$ is described by its
{\em code} $c_P\in\{0,1\}^n$, which is given by the first $n$ symbols
of the itinerary $k(p)$ of its rightmost point~$p$: thus, for example,
the period 5 orbit which contains the point with itinerary
$\overline{01001}$ has code $10010$. A word $w\in\{0,1\}^n$ is
therefore the code of a period $n$ horseshoe orbit if and only if it
is {\em maximal}: that is, the infinite repetition
$\overline{w}\in\{0,1\}^\Nset$ of $w$ is strictly greater than its
shifts $\sigma^i(\overline{w})$ in the unimodal order for $1\le i<n$.

The following definitions are due to Boyland~\cite{Boy,Boy2}. Let
$f\colon D^2\to D^2$ and $g\colon D^2\to D^2$ be
orientation-preserving homeomorphisms having periodic orbits $P$ and
$Q$ respectively. Then $(P,f)$ and $(Q,g)$ {\em have the same braid
type} if there is a homeomorphism $h\colon (D^2,P)\to(D^2,Q)$ such
that $f\colon(D^2,P)\to(D^2,P)$ is isotopic (rel. $P$) to $h\I\circ
g\circ h\colon(D^2,P)\to(D^2,P)$ (if either $P$ or $Q$ lies on
$\partial D^2$, then the corresponding homeomorphism should first be
extended arbitrarily over an exterior collar). The {\em braid type}
$\bt(P,f)$ of $(P,f)$ is its equivalence class under this
relation. Since the braid types of period $n$ orbits correspond to
conjugacy classes in the mapping class group of the $n$-punctured
disk, they can be classified as finite order, reducible, or
pseudo-Anosov by means of Thurston's
classification~\cite{Thu}. Boyland's {\em forcing relation} $\le$ on
the set $\BT$ of braid types is defined as follows: $\beta\le\gamma$
if and only if every orientation-preserving homeomorphism $f\colon
D^2\to D^2$ which has a periodic orbit of braid type $\gamma$ also has
one of braid type $\beta$. Boyland proved~\cite{Boy,Boy2} that $\le$
is a partial order on $\BT$.

An alternative characterisation of braid type indicates the dynamical
significance of the definition: $(P,f)$ and $(Q,g)$ have the same
braid type if and only if there exists an isotopy $\{f_t\}$ from $f$
to $g$ and a path $P_t$ in $(D^2)^n$ from $P$ to~$Q$ such that $P_t$
is a (least) period $n$ orbit of $f_t$ for all $t$.

The braid type $\bt(P,F)$ of a horseshoe periodic orbit $P$ will here
be denoted simply by $\bt(P)$. It is well known that two periodic
orbits $P$ and $Q$ of the horseshoe whose codes $c_P$ and $c_Q$ differ
only in their final symbol have the same braid type. Thus, for
example, the two orbits with codes $10010$ and $10011$ have the same
braid type; the code of either one of these orbits is often written
$c_P=1001\opti$ to reflect the fact that the distinction between the
two is unimportant in so far as braid type is concerned.

The conjecture presented in~\cite{dCH2} is based upon the parsing of
the code of any horseshoe periodic orbit $P$ which is not of finite
order braid type into two parts: the {\em prefix} and the {\em
decoration}. In order to define these, it is necessary first to
describe the {\em height} $q(P)\in(0,1/2]\cap\Qset$ of $P$: this is an
invariant of braid type which plays a central role in the
conjecture. Motivation for the definition is given in~\cite{Ha}, and a
program for computing heights of horseshoe periodic orbits can be
found at~\cite{programs}.

\begin{defns}
Let $c\in\{0,1\}^\Nset$ be a sequence containing infinitely many
$1$s. If $c$ starts with $0$ or with $11$, then define the {\em height}
$q(c)$ of $c$ to be $1/2$. Otherwise, write
$$c=10^{\kappa_1}1^{\mu_1}0^{\kappa_2}1^{\mu_2}\ldots,$$
where each $\kappa_i\ge0$, each $\mu_i$ is either $1$ or $2$, and
$\mu_i=1$ only if $\kappa_{i+1}>0$ (thus $\kappa_i$ and $\mu_i$ are
uniquely determined by $c$). For each $r\ge 1$, define

$$I_r(c)=\left(\frac{r}{2r+\sum_{i=1}^r\kappa_i},
\frac{r}{(2r-1)+\sum_{i=1}^r\kappa_i}\right].$$
If there is a (least) integer $s\ge 1$  such that either $\mu_s=1$ or
$\bigcap_{i=1}^{s+1}I_i(c)=\emptyset$, then let
$\bigcap_{i=1}^sI_i(c)=(x,y]$. The height $q(c)$ of $c$ is given by
\[q(c)=\left\{
\begin{array}{ll}
x & \mbox{ if }\mu_s=2\mbox{ and }w<z\mbox{ for all }w\in I_{s+1}(c)
\mbox{ and }z\in \bigcap_{i=1}^sI_i(c)\\
y & \mbox{ if }\mu_s=1\mbox{, or } \mu_s=2\mbox{ and }
w>z\mbox{ for all }w\in I_{s+1}(c)\mbox{ and }z\in\bigcap_{i=1}^s I_i(c).
\end{array}
\right.\]

Let $P$ be a horseshoe periodic orbit with code $c_P$. If the
sequence $\overline{c_P}\in\{0,1\}^\Nset$ does not contain the word
$010$, then change the final symbol of $c_P$ from $1$ to $0$. Then the
{\em height} $q(P)$ of $P$ is defined to be $q(\overline{c_P})$.
\end{defns}

Notice that, having changed the final symbol of $c_P$ if necessary,
some $\mu_i$ is equal to $1$ and hence the above algorithm
terminates. The definition of height can be extended to elements of
$\{0,1\}^\Nset$ which contain only finitely many $1$s, or for which
the above algorithm does not terminate (see for example
Lemma~\ref{lem:riprops}~c)), although such cases will not arise in
this paper. With this extension, the function
$q\colon\{0,1\}^\Nset\to[0,1/2]$ is decreasing (but not strictly
decreasing) with respect to the unimodal order on $\{0,1\}^\Nset$ and
the usual order on $[0,1/2]$.

\begin{exmp}
\label{ex:height}
Let $P$ be the period $17$ orbit with code $10011011001011010$. Then
$\kappa_1=2$, $\mu_1=2$, $\kappa_2=1$, $\mu_2=2$, $\kappa_3=2$, and
$\mu_3=1$. Thus $I_1=(1/4,1/3]$, $I_2=(2/7,2/6]$, and
$I_3=(3/11,3/10]$. Since $\mu_3=1$, the algorithm terminates with
$\bigcap_{i=1}^3I_i=(2/7,3/10]$, and hence $q(P)=3/10$.
\end{exmp}

\begin{figure}
\lab{0}{0}{}
\lab{1}{1}{}
\lab{2}{(0,0)}{r}
\lab{3}{(10,3)}{l}
\begin{center}
\pichere{0.4}{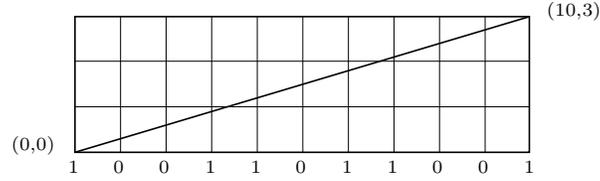}
\end{center}
\caption{$c_{3/10}=10011011001$}
\label{fig:linebox}
\end{figure}

The following definition and theorem are taken from~\cite{Ha}. Given
$q=m/n\in(0,1/2]$, define $c_q\in\{0,1\}^{n+1}$ by
$$c_q=10^{\kappa_1}1^20^{\kappa_2}1^2\ldots 1^20^{\kappa_m}1,$$
where
\begin{equation}
\label{eq:kappa}
\kappa_i=\left\{
\begin{array}{ll}
\lfloor n/m\rfloor -1 & \mbox{ if }i=1\\
\lfloor in/m\rfloor - \lfloor (i-1)n/m\rfloor - 2 & \mbox{ if }2\le
i\le m
\end{array}
\right.
\end{equation}
(here $\lfloor x\rfloor$ denotes the greatest integer which does not
exceed $x$). An alternative description of the words $c_q$ which is
sometimes useful is as follows: consider the straight line $L$ in
$\Rset^2$ from $(0,0)$ to $(n,m)$. For $0\le i\le n$, let $s_i$ be
equal to $1$ if $L$ crosses a horizontal line $y=\text{integer}$ for
$i-1<x<i+1$, and equal to $0$ otherwise. Then $c_q=s_0s_1\ldots
s_{n}$. Thus the example of Fig.~\ref{fig:linebox} shows that
$c_{3/10}=10011011001$.

It is immediate from this description that the words $c_q$ are
palindromic: that is, $\kappa_i=\kappa_{m+1-i}$ for all $i$.

\begin{thm}
\label{thm:hs}
\begin{enumerate}[a)]
\item Let $P$ be a horseshoe periodic orbit with height $q(P)=m/n$ in lowest
terms. Then $P$ has period $n$ if and only if it has finite order
braid type: in this case, $F$ is isotopic rel.~$P$ to a rigid rotation
through $2\pi m/n$.  Otherwise, the period of $P$ is at least $n+2$,
and $c_P$ starts with the word~$c_{m/n}$.
\item Let $P$ and $R$ be horseshoe periodic orbits. If
$\bt(P)\ge\bt(R)$ then $q(P)\le q(R)$. In particular, height is an
invariant of braid type.
\end{enumerate}
\end{thm}

Thus the code $c_P$ of a horseshoe periodic orbit which is not of
finite order braid type can be written $c_P=c_{q(P)}v$ for some word
$v$ of length at least $1$. This makes possible the following
definitions, which are taken from~\cite{dCH2}:
\begin{defns}
Let $P$ be a period $N$ orbit of the horseshoe which is not of finite
order braid type, with height $q=q(P)=m/n$. The {\em prefix} of $P$ is
the word $c_q$. The {\em decoration} of $P$ is defined to be~$\ast$ if
$N=n+2$, and to be the element $w$ of $\{0,1\}^{N-n-3}$ such that
\[c_P=c_q\opti w\opti\] otherwise.
\end{defns}

\begin{exmp} 
Let $P$ be the period $17$ orbit with code
$c_P=10011011001011010$. Then $q(P)=3/10$ as shown in
Example~\ref{ex:height}. Hence $P$ has prefix $10011011001=c_{3/10}$,
and decoration $1101$.
\end{exmp}

Only certain heights $q$ are compatible with a given decoration $w$
(namely those for which $c_q\opti w\opti$ is a maximal word, and hence
describes the rightmost point of a periodic orbit). The following
lemma (from~\cite{dCH2}) gives the compatibility conditions:
\begin{lem}
\label{lem:qw}
Let $w$ be a decoration, and define $q_w\in(0,1/2]\cap\Qset$ by
$q_\ast=1/2$ and
\[q_w=\min_{0\le i\le k+2}q\left(\sigma^i\left(\overline{10w0}\right)
\right)\] if $w\in\{0,1\}^k$.
Then each of the four words $c_q\opti w\opti$ (or each of the two
words $c_q\opti$ when $w=\ast$) is maximal of height $q$ when
$0<q<q_w$, and none is maximal of height $q$ when $q_w<q\le 1/2$.
\end{lem}

The reason that the four orbits with codes $c_q\opti w\opti$ are
considered together is that they all have the same braid type: the
following general result will appear in~\cite{invs}. That the braid
type is unchanged on changing the final symbol of the code is a
triviality: the content of the theorem is that changing the symbol
between prefix and decoration also leaves the braid type unchanged.

\begin{thm}
\label{thm:general}
Let $w\not=\ast$ be a decoration, and $q\in(0,q_w)\cap\Qset$. Then the
four horseshoe periodic orbits with codes $c_q\opti w\opti$ have the
same braid type.
\end{thm}

The following conjecture is motivated and stated in~\cite{dCH2}. It
involves the notion of two horseshoe periodic orbits of pseudo-Anosov
braid type having the {\em same topological train track type}. A
formal definition of this term can be found in~\cite{dCH2}: however
its intuitive meaning should be clear in the context of this paper (in
which `star' train track types are considered).
\begin{conj}
\label{conj:c1}
\begin{enumerate}[a)]
\item Let $P$ be a horseshoe periodic orbit with height $q$ and
decoration $w$. If $q\not=q_w$ then $P$ has pseudo-Anosov braid type.
\item Let $P'$ be another periodic orbit with height $q'$ and
decoration $w$. Then $\bt(P')\le\bt(P)$ if and only if $q\le q'$
(i.e. the set of braid types of periodic orbits with decoration $w$ is
totally ordered by the forcing relation). Moreover, if $q'<q_w$ then
$P$ and $P'$ have the same topological train track type.
\item There is an equivalence relation $\sim$ on the set $\cD$ of
decorations with the property that two horseshoe periodic orbits have
the same braid type if and only if they have equal heights and
equivalent decorations.
\item There is a partial order $\preceq$ on $\cD/\!\!\sim$ with the
property that if $P$ and $P'$ are periodic orbits with heights $q$,
$q'$ and decorations $w$, $w'$ such that $q<q'$ and $[w']\preceq[w]$,
then $\bt(P')\le \bt(P)$.
\end{enumerate}
\end{conj}

Notice that part~b) trivializes the problem of comparing braid types
within a family of fixed decoration: the forcing order is given by the
unimodal order on the codes of the periodic orbits.  As explained
in~\cite{dCH2}, parts~a) and~b) of this conjecture (i.e., the
statements which only concern a single decoration) can be proved for
certain particular choices of decoration~$w$. In this paper an
infinite family of decorations is considered, making it possible to
address parts~c) and~d) of the conjecture meaningfully. The
decorations concerned are those which give rise to `star' topological
train track types, and are given by words $w_q$ with
$q\in(0,1/2]\cap\Qset$ defined as follows:

\begin{defn}
Let $q=m/n\in(0,1/2]\cap\Qset$. If $q\not=1/2$ then $w_q$ is the element of
$\{0,1\}^{n-3}$ obtained by deleting the initial symbols $10$ and the
final symbols $01$ from $c_q$. If $q=1/2$ then $w_q=\ast$.
\end{defn}

Thus, for example, $w_{1/3}$ is the empty decoration, $w_{1/4}=0$,
$w_{1/5}=00$, and $w_{2/5}=11$. These decorations satisfy $q_{w_q}=q$,
are mutually non-equivalent under $\sim$, and are totally ordered by
$\preceq$, with $w_q\preceq w_{q'}$ if and only if $q\le q'$. Thus the
proof of Conjecture~\ref{conj:c1} above for these decorations gives
the following theorem, which summarizes some of the main results of
this paper.

\begin{thm}
\label{thm:summ}
Let $r,r'\in(0,1/2]\cap\Qset$ and $q,q'\in(0,1/2)\cap\Qset$ with $q\le r$
and $q'\le r'$. Let $P$ be a horseshoe periodic orbit of height $q$
and decoration $w_r$, and $P'$ be a horseshoe periodic orbit of height
$q'$ and decoration $w_{r'}$. Then
\begin{enumerate}[a)]
\item If $q\not=r$ then $P$ has pseudo-Anosov braid type.
\item If $r=r'$ and $q,q'<r$ then $P$ and $P'$ have the same (star)
topological train track type. Moreover $\bt(P')\le\bt(P)$ if and only
if $q\le q'$.
\item $P$ and $P'$ have the same braid type if and only if $q=q'$ and
$r=r'$.
\item If $q<q'$ and $r\ge r'$ then $\bt(P')\le\bt(P)$.
\end{enumerate}
\end{thm}

Part~a) of this theorem is given by Corollary~\ref{cor:pA} below,
Part~b) by Theorem~\ref{thm:main}, and Part~c) by
Theorem~\ref{thm:general} and Lemma~\ref{lem:identify}, as described
following the statement of Lemma~\ref{lem:identify}.  Part~d) is
proved using quite different techniques from those of this paper: a
proof will appear in~\cite{invs}.

Lemma~\ref{lem:identify} below will be used to identify periodic
orbits with decoration $w_{m/n}$ for some $m/n\in(0,1/2]$. Its proof
is relatively long and technical, and uses methods from~\cite{Ha}
which are not required elsewhere in this paper --- it has therefore
been relegated to Section~\ref{sec:app}.

\begin{defns}
Let $P$ be a horseshoe periodic orbit of period $N>1$. The {\em
rotation number} $\rho(P)\in(0,1/2]\cap\Qset$ of $P$ is its
$F$-rotation number about the fixed point with code $1$. The {\em
rotation interval} of $P$ is the set
\[\rho i(P)=\{\rho(P')\colon \bt(P')\le\bt(P)\}.\]
\end{defns}

The rotation number of a horseshoe periodic orbit is an invariant of
braid type~\cite{Ha} (this statement is not quite as obvious as it may
at first appear), and hence so also is the rotation interval.  By a
theorem of Handel~\cite{Handelrot}, $\rho i(P)$ is equal either to
$\{\rho(P)\}$ or to a set of the form $[a,b]\cap\Qset$ (which will
here be denoted simply $[a,b]$) for some $a<b\in\Qset$.

\begin{lem}
\label{lem:identify}
Let $P$ be a period $N>1$ orbit of the horseshoe with non-trivial
rotation interval $\rho i(P)=[u/v,m/n]$. Then $N\ge v+n$. Moreover, 
$N=v+n$ if and only if $c_P$ is one of the four words $c_{u/v}\opti
w_{m/n}\opti$ (or one of the two words $c_{u/v}\opti$ in the case
$m/n=1/2$). 
\end{lem}

In other words, if the period of $P$ is equal to the sum of the
denominators of the endpoints of its rotation interval, then $P$ has
height $u/v$ and decoration $w_{m/n}$. Since the rotation interval is
a braid type invariant, it follows that if $P$ has height $q$ and
decoration $w_r$, and $P'$ is a horseshoe periodic orbit with the same
braid type as $P$, then $P'$ also has height $q$ and decoration
$w_r$. In particular, this establishes (and extends) the `only if'
part of Theorem~\ref{thm:summ}~c). The `if' part follows from
Theorem~\ref{thm:general}. Moreover, the following corollary provides
a proof of Theorem~\ref{thm:summ}~a).

\begin{cor}
\label{cor:pA}
Let $u/v<m/n\le 1/2$. Then the horseshoe periodic orbits of height
$u/v$ and decoration $w_{m/n}$ have pseudo-Anosov braid type.
\end{cor}
\begin{proof}
Let $P$ be a periodic orbit of height $u/v$ and decoration
$w_{m/n}$. Then $\rho i(P)=[u/v,m/n]$. Periodic orbits of finite order
type have trivial rotation intervals. If $P$ had reducible braid type,
then a horseshoe periodic orbit of period less than $v+n$ (namely one
given by the isotopy class corresponding to the outermost reducible
component of $F\colon (D^2,P)\to(D^2,P)$) would also have rotation
interval $[u/v,m/n]$, contradicting Lemma~\ref{lem:identify}.
\end{proof}

Lemma~\ref{lem:identify} may be of independent interest. The authors
know of no pseudo-Anosov homeomorphism $f$ of the annulus, relative to
a single periodic orbit~$P$, for which the sum of the denominators of
the endpoints of the rotation interval of $f$ is less than the period
$n$ of $P$. If this more general result could be proved, then
(together with the fact~\cite{BGH} that the rotation number $m/n$ of
$P$ is contained in the interior of the rotation interval of $f$) it
would provide a natural generalization of Boyland's
theorem~\cite{Boyri}, which gives bounds on the size of the rotation
interval in the case where $m$ and $n$ are coprime.

\subsection{The Farey Graph}
\label{sec:farey}
In this section some well known results about Farey sequences and the
Farey graph are presented, some notation is introduced, and some
simple number-theoretic lemmas which will be used later are stated. In
the usual treatment (see for example~\cite{Farey}, where the
assertions made in this section are proved), the Farey sequences
contain rationals between $0$ and $1$: since only rationals between
$0$ and $1/2$ are of interest in this paper, the definitions have been
modified accordingly. Throughout the paper, all rationals are assumed
to be written in lowest terms.

\begin{defns}
Let $n\ge1$ be an integer. The {\em Farey sequence $\cF_n$ of order
$n$} is the (finite) sequence of rational numbers in $[0,1/2]$ whose
denominators do not exceed $n$, arranged in ascending order. Two
rationals $h/k$ and $h'/k'$ in $[0,1/2]$ are {\em Farey neighbours} if
there is some $n$ such that they are consecutive elements of
$\cF_n$. Given a rational $h/k\in(0,1/2)$, the {\em left and right
Farey parents} $\LFP(h/k)$ and $\RFP(h/k)$ are the elements of $\cF_k$
which precede and follow $h/k$; and $h/k$ is said to be the {\em Farey
child} of its parents.
\end{defns}
In particular, both $\LFP(h/k)$ and $\RFP(h/k)$ are Farey neighbours
of $h/k$. It is well known that if $h/k<h'/k'$ are Farey neighbours, then
$h'k-hk'=1$; and that if $\LFP(h/k)=u/v$ and $\RFP(h/k)=p/q$, then
$h=u+p$ and $k=v+q$.

\begin{figure}[htbp]
\lab{1}{\frac{0}{1}}{t}
\lab{2}{\frac{1}{3}}{t}
\lab{3}{\frac{1}{2}}{t}
\lab{4}{\frac{1}{4}}{t}
\lab{5}{\frac{2}{5}}{t}
\lab{6}{\frac{1}{5}}{t}
\lab{7}{\frac{1}{6}}{t}
\lab{8}{\frac{2}{7}}{t}
\lab{9}{\frac{2}{9}}{t}
\lab{a}{\frac{3}{11}}{t}
\lab{b}{\frac{3}{10}}{t}
\lab{c}{\frac{3}{8}}{t}
\lab{d}{\frac{3}{7}}{t}
\lab{e}{\frac{4}{11}}{t}
\lab{f}{\frac{5}{13}}{t}
\lab{g}{\frac{5}{12}}{t}
\lab{h}{\frac{4}{9}}{t}
\begin{center}
\pichere{0.6}{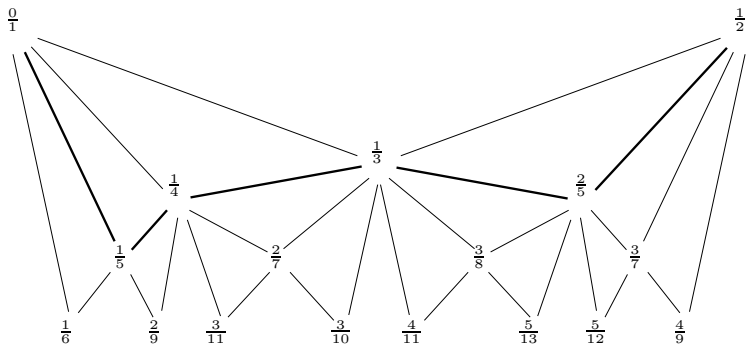}
\end{center}
\caption{Part of the Farey Graph}
\label{fig:fareygra}
\end{figure}

The rationals in $[0,1/2]$ can be organized as the vertices of the
{\em Farey Graph}, in which an edge joins two vertices if and only if
one of the associated rationals is the left or right Farey parent of
the other. Part of the Farey graph is depicted in
Fig.~\ref{fig:fareygra}. Notice that every vertex has infinite
valence (that is, every rational is the parent of infinitely many
other rationals). The {\em immediate left (respectively right) Farey
child} of a rational $m/n$ is the rational $h/k$ of smallest
denominator for which $\LFP(h/k)=m/n$ (respectively $\RFP(h/k)=m/n$):
it is the child of $m/n$ and $\LFP(m/n)$ (respectively
$\RFP(m/n)$). Thus, for example, the immediate children of $2/5$ are
$3/8$ and $3/7$. The Farey sequences appear as finite connected
subtrees of the Farey graph: for example, $\cF_5$ is shown in the
figure with bolder lines.

\begin{defn}
Let $m/n$ be a rational in $(0,1/2)$. Then the {\em left Farey
sequence} $\LFS(m/n)$ of $m/n$ is the (finite) sequence
$(0=u_1/v_1,u_2/v_2,\ldots,u_\alpha/v_\alpha)$, where
$u_\alpha/v_\alpha=LFP(m/n)$, and $u_i/v_i=LFP(u_{i+1}/v_{i+1})$ for
$1\le i< \alpha$.
\end{defn}

Thus, for example, $\LFS(3/10)=(0,1/4,2/7)$.

\begin{defn}
Let $m/n$ be a rational in $(0,1/2)$, and denote addition modulo $n$
by $+_n$. Given $r,s\in\Zset_n$, let
\[\cO_{m/n}[r,s]=\{r+_njm\colon\,0\le j\le k,\text{ where $k\ge0$ is
least such that $r+_nkm=s$}\}.\]
\end{defn}
That is, $\cO_{m/n}[r,s]$ is the shortest segment of the orbit of $r$
in $\Zset_n$ under addition of $m$ which ends with $s$. When $r\le s$ are
integers, the notation $[r,s]$ will also be used as a shorthand for
$\{r,r+1,\ldots,s\}$. 

The following simple Lemma contains the results about the Farey
graph which will be needed in the remainder of the paper. 
\begin{lem}
\label{lem:fareyprops}
Let $m/n\in(0,1/2)$, and suppose that $\LFP(m/n)=u/v$ and
$\RFP(m/n)=p/q$. Then
\begin{enumerate}[a)]
\item $\card \cO_{m/n}[m,n-1]=q=n-v$, and
$\card\cO_{m/n}[m-1,0]=n-q=v$.
\item
$\card(\cO_{m/n}[m,n-1]\cap[1,m])=\card(\cO_{m/n}[m,n-1]\cap[n-m,n-1])=p=m-u$.
\item $\card(\cO_{m/n}[m-1,0]\cap[1,m])=
\card(\cO_{m/n}[m-1,0]\cap[n-m,n-1]=m-p=u$.
\end{enumerate}
\end{lem}
\begin{proof}
\begin{enumerate}[a)]
\item By definition of $R=\cO_{m/n}[m,n-1]$ its cardinality is $k+1$,
where $k\in\Zset_n$ is such that $m+km\equiv n-1\bmod{n}$: that is,
$(k+1)(m-n)\equiv -1\bmod{n}$, so $k+1$ is the multiplicative inverse
$(n-m)_n\I$ of $n-m$ in $\Zset_n$.

Since $p/q=\RFP(m/n)$, it follows that $q\in[1,n-1]$, and that
$qm=np-1$. Thus $q(n-m)\equiv -qm\equiv 1\bmod{n}$, so $(n-m)_n\I=q$
as required. However $v+q=n$, so $\card R=n-v$ also. Since $R$ and
$S=\cO_{m/n}[m-1,0]$ partition $\Zset_n$, the cardinality of $S$
follows.

\item Let $\card\left(R\cap[1,m]\right)=s$. Since 
$R$ has cardinality $q$, it follows that $m+(q-1)m=n-1+(s-1)n$, or
$qm+1=sn$. Thus $s=p$. However $u+p=m$, so $s=m-u$ also. A similar
argument applies to $R\cap[n-m,n-1]$.

\item Immediate from part~b). 
\end{enumerate}
\end{proof}

The function defined in the following lemma will play a central
role in this paper.
\begin{lem}
\label{lem:xi}
Let $u/v<p/q$ be Farey neighbours. Then the function
\[\xi_{u/v,p/q}\colon(0,1)\cap\Qset\to(u/v,p/q)\cap\Qset\] defined by
\[\xi_{u/v,p/q}\left(\frac{r}{s}\right)=\frac{rp+(s-r)u}{rq+(s-r)v}\]
is an increasing bijection.
\end{lem}
\begin{proof}
Since $u/v<p/q$ are Farey neighbours, $pv-qu=1$. Hence if
$m/n\in(u/v,p/q)$, then $m/n=\xi_{u/v,p/q}((vm-un)/((v-q)m+(p-u)n))$
(note that $0<vm-un<(v-q)m+(p-u)n=(vm-un)+(pn-qm)$ (since
$u/v<m/n<p/q$), and that $vm-un$ and $(v-q)m+(p-u)n$ are coprime
(since $m=(vm-un)(p-u)+((v-q)m+(p-u)n)u$ and
$n=(vm-un)(q-v)+((v-q)m+(p-u)n)v$ are coprime)). Thus $\xi_{u/v,p/q}$
is surjective. It is also strictly increasing, since if $0<r/s<a/b<1$
then
$\xi_{u/v,p/q}(a/b)-\xi_{u/v,p/q}(r/s)=\frac{as-br}{(aq+(b-a)v)(rq+(s-r)v)}>0$.
\end{proof} 
Notice in particular that $\xi_{u/v,p/q}(1/2)=(u+p)/(v+q)$ is the
Farey child of $u/v$ and $p/q$.

\section{Periodic orbits of star maps}
\label{sec:perorb}
\subsection{Star maps and $\ast$-orbits}

\begin{figure}[htbp]
\lab{1}{e_0}{t}
\lab{2}{e_1}{tl}
\lab{3}{e_2}{bl}
\lab{4}{e_3}{b}
\lab{5}{e_4}{r}
\lab{6}{f_{2/5}(e_3)}{b}
\lab{7}{f_{2/5}(e_4)}{br}
\lab{8}{f_{2/5}(e_0)}{t}
\lab{9}{f_{2/5}(e_1)}{t}
\lab{0}{f_{2/5}(e_2)}{l}
\lab{a}{f_{2/5}}{b}
\begin{center}
\pichere{0.6}{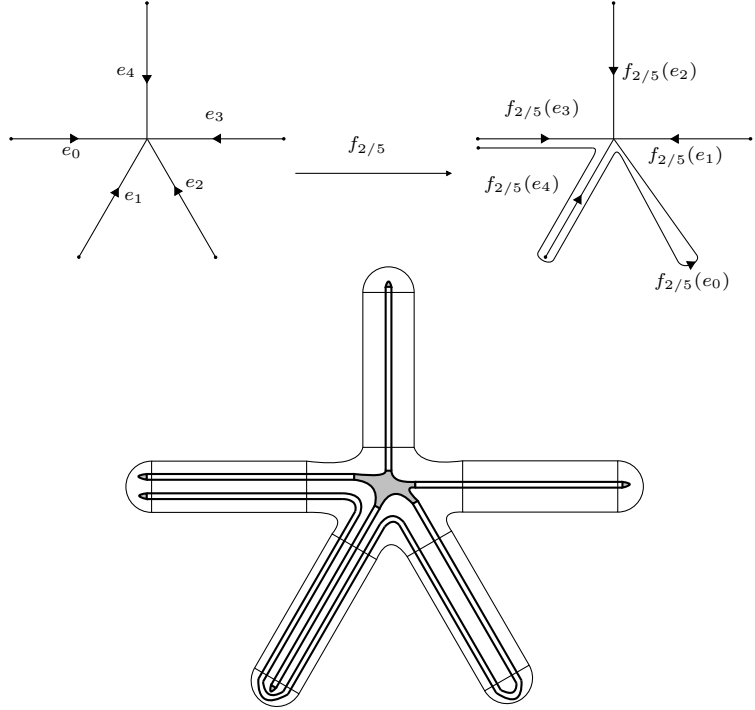}
\end{center}
\caption{The tree map $f_{2/5}:\Gamma_5\to\Gamma_5$ 
and the thick tree map $F_{2/5}:T_5\to T_5$}
\label{fig:thick}
\end{figure}

For each $n\ge 2$ let $\Gamma_n\subseteq D^2$ be an {\em $n$-star}:
that is, a tree with $n$ edges $e_0,\ldots,e_{n-1}$ of equal length,
each of which has a valence $1$ vertex at its initial point, and which
meet at a valence $n$ vertex $v$ at their final points, cyclically
ordered according to their indices. For each rational $m/n\in(0,1/2]$
(written in its lowest terms), let $f_{m/n}\colon\Gamma_n\to\Gamma_n$
be the tree map (see Fig.~\ref{fig:thick}) with image edge paths
\begin{eqnarray*}
f_{m/n}(e_0)&=&e_0\be_1 e_1\be_2e_2 \ldots \be_m e_m\\
f_{m/n}(e_r)&=&e_{r+_nm}\qquad\mbox{($r\not=0$)},
\end{eqnarray*}
which expands each edge uniformly away from the preimages of vertices
(each connected component of which is either a nontrivial interval or,
in the case of the preimage component containing $v$, is homeomorphic
to $\Gamma_n$). Figures depicting $f_{m/n}\colon\Gamma_n\to\Gamma_n$
(with the exception of Fig.~\ref{fig:scheme}) are always drawn with
$e_0$ and $e_{n-m}$ horizontal.

Let $F_{m/n}\colon T_n\to T_n$ be a {\em thick tree map} (see for
example~\cite{dCH1} for a formal definition) corresponding to
$f_{m/n}\colon
\Gamma_n\to\Gamma_n$ (see Fig.~\ref{fig:thick}). 
Thus $T_n\subseteq D^2$ is a topological disk, and there is a
continuous map $p\colon T_n\to\Gamma_n$ such that $p\I(x)$ is a disk
if $x$ is a vertex of $\Gamma_n$, and an interval otherwise. The map
$F_{m/n}$ is an embedding which contracts each such {\em decomposition
element} $p\I(x)$ into a decomposition element, in such a way that it
induces $f_{m/n}$ on $\Gamma_n$ (so $F_{m/n}(p\I(x))\subset
p\I(f_{m/n}(x))$ for all $x\in\Gamma_n$). Where necessary, $F_{m/n}$
is considered as a homeomorphism $D^2\to D^2$, by extending from $T_n$
without introducing any new periodic orbits.

There is a natural one-to-one correspondence (which will be invoked
without comment in the remainder of the paper) between the periodic
orbits of $f_{m/n}$ and those of $F_{m/n}$ (with a periodic point $x$
of $f_{m/n}$ corresponding to the unique periodic point of $F_{m/n}$
in $p\I(x)$).  Notice that $F_{1/2}$ is the horseshoe map $F$ after
blowing up the leaf containing the fixed point of code $1$ into a
disk: in particular, there is a braid type preserving bijection
between the set of periodic orbits of $F_{1/2}$ and the set of
periodic orbits of $F$.

In this paper only periodic orbits $P$ which satisfy
certain non-triviality conditions are considered:

\begin{defn}
A periodic orbit $P$ of $f_{m/n}\colon \Gamma_n\to\Gamma_n$ is called
a {\em $\ast$-orbit} if
\begin{enumerate}[$\ast$\,a)]
\item $P\not=\{v\}$.
\item $P\cap e_r\not=\emptyset$ for all $r$.
\item $f_{m/n}(P\cap e_0)\not\subseteq e_m$.
\item If $p_r$ denotes the point of $P\cap e_r$ closest to the initial
point of $e_r$, then $f_{m/n}(p_r)=p_{r+_nm}$ for all $r\not=0$.
\end{enumerate}
\end{defn}

The reasons for imposing three of these conditions are intuitively
clear: a)~states that $P$ is not a fixed point; b)~that it explores
each of the edges of $\Gamma_n$; and c)~that the points of $P$ in
$e_0$ have images in more than one edge: if this were not true, the
braid type of $P$ would either be finite order (if each edge contained
just one point of $P$) or reducible (with $n$ reducing curves, each
bounding a disk containing the points of $P$ on one of the edges of
$\Gamma_n$). Condition~d) is less clear: the motivation is that when
$f_{m/n}$ is `truncated' with respect to~$P$ (see
Definition~\ref{def:trunc}), there is only one point of $P$ (namely the
preimage of~$p_m$) at which $f_{m/n}$ is not locally
injective. Lemma~\ref{lem:span} below is one important consequence of
this condition.

\begin{defn}
Let $P$ be a periodic orbit of $f_{m/n}$. Then the {\em span}
$\Gamma_n^P$ of $P$ is the smallest connected subset of $\Gamma_n$
containing $P$.
\end{defn}

\begin{lem}
\label{lem:span}
Let $m/n\in(0,1/2]$. Then the set of spans of $\ast$-orbits of
$f_{m/n}$ is totally ordered by inclusion.
\end{lem}
\begin{proof}
Let $P$ and $P'$ be $\ast$-orbits of $f_{m/n}$, and for each $0\le
r<n$ let $p_r,p_r'\in e_r$ be the points of the orbits closest to the
initial point of $e_r$ (which exist by $\ast$\,b)). If $p_m$ is closer
to the initial point of $e_m$ than $p'_m$, then applying $\ast$\,d)
inductively gives that $p_r$ is closer than $p'_r$ to the initial
point of $e_r$ for all~$r$, and hence
$\Gamma_n^{P'}\subseteq\Gamma_n^P$.
\end{proof}

The following result, which is contained in Theorem~\ref{thm:hscode}
below, makes it possible to identify the braid types of $\ast$-orbits
using Lemma~\ref{lem:identify}. It reflects the fact that the thick
tree maps $F_{m/n}\colon T_n\to T_n$ can be obtained from the
horseshoe by pruning (i.e. by destroying some dynamics), as described
in Section~\ref{sec:pfh}.

\begin{lem}
\label{lem:starhs}
Every $\ast$-orbit of $f_{m/n}$ has the braid type of some periodic
orbit of the horseshoe.
\end{lem}

\subsection{Describing $\ast$-orbits}
In this section a combinatorial method for describing $\ast$-orbits is
developed. Since $f_{m/n}(e_r)=e_{r+_nm}$ for $r\not=0$, the main work
required is in describing the images of the points of the orbit in
$e_0$. The first step is to make precise the notion of truncating
$f_{m/n}$ with respect to a $\ast$-orbit $P$.

\begin{defn}
\label{def:trunc}
Let $P$ be a $\ast$-orbit of $f_{m/n}$. Let
$r_P\colon\Gamma_n\to\Gamma_n^P$ be the map defined by $r_P(x)=x$ if
$x\in\Gamma_n^P$, and $r_P(x)$ is the endpoint of $\Gamma_n^P$
contained on the same edge of $\Gamma_n$ as $x$ if
$x\not\in\Gamma_n^P$. The {\em truncation} $f_{m/n}^P$ of $f_{m/n}$
with respect to $P$ is the map $f_{m/n}^P=r_P\circ f_{m/n}\colon
\Gamma_n^P\to\Gamma_n^P$.
\end{defn}

The following definitions give a basic classification of $\ast$-orbits.

\begin{defns}
\label{def:abc}
Let $P$ be a $\ast$-orbit of $f_{m/n}$ with $\card(P\cap e_0)=N$, and
label the points of $P\cap e_0$ as $p_0,p_1,\ldots,p_{N-1}$ in order
from the initial to the final point of $e_0$. Define a partition
\[\{0,\ldots,N-1\}=A\cup B\cup C\] by
\[i\in\left\{\begin{array}{ll}
A &\mbox{ if $f^P_{m/n}$ is locally orientation-reversing at $p_i$}\\
B &\mbox{ if $f^P_{m/n}$ is locally orientation-preserving at $p_i$}\\
C &\mbox{ if $f^P_{m/n}$ is not locally injective at $p_i$.}
\end{array}\right.\]

Observe that $i\in C$ if and only if $f_{m/n}(p_i)$ is the point of $P$
closest to the initial point of $e_m$, and hence $\card C=1$.

Given a rational $m/n\in(0,1/2]$, an integer $k\in\{0,\ldots,m-1\}$,
and $\gamma\in\{A,B\}$ (with $\gamma=B$ if $k=0$), write
$\cP(m/n,k,\gamma)$ for the set of all $\ast$-orbits of $f_{m/n}$ with
$f(p_0)\in e_k$ and $0\in\gamma$,
$\cP(m/n,k)=\cP(m/n,k,A)\cup\cP(m/n,k,B)$, and
$\cP(m/n)=\bigcup_k\cP(m/n,k)$, the set of all $\ast$-orbits of $f_{m/n}$.
\end{defns}

The rest of the description of a $\ast$-orbit is contained in the next
definitions. 

\begin{defns}
For each integer $n\ge 2$, define $\cD_n$ to be the set of all
triples
\[d=\left((N_0,\ldots,N_{n-1}),\pi,(A,B,C)\right), \]
where $N_r$ is a positive integer for $0\le r<n$; $\pi$ is a cyclic
permutation of the set
\[\cL=\cL_d=\{(r,s)\,:\,0\le r<n,\,0\le s<N_r\};\]
and $(A,B,C)$ is a partition of $\{0,\ldots,N_0-1\}$ with
$\card C=1$.

Let $P\in\cP(m/n,k,\gamma)$. Then the {\em data} of $P$ is the element
\[d(P)=\left((N_0,\ldots,N_{n-1}),\pi,(A,B,C)\right)\]
of $\cD_n$ obtained as follows: $N_r=\card(P\cap e_r)$ for each $r$
with $0\le r<n$. Identify $P$ with $\cL=\cL_{d(P)}=\{(r,s)\,:\,0\le
r<n,\,0\le s<N_r\}$ by labelling the points of $P\cap e_r$ as $(r,0),
(r,1),\ldots(r,N_r-1)$ from the initial to the final point of $e_r$,
and let $\pi=f|_P\colon\cL\to\cL$. Finally, let $(A,B,C)$ be the
partition of $\{0,\ldots,N-1\}$ given by Definitions~\ref{def:abc}.
\end{defns}

Having the same data is the basic equivalence relation which says that
two \mbox{$\ast$-orbits} have the same `shape'. In particular, it is
clear that two $\ast$-orbits with the same data have the same braid
type as periodic orbits of $F_{m/n}\colon T_n\to T_n$. In this paper
no distinction is made between orbits with the same data: thus, for
example, the statement that $f_{m/n}$ has exactly one $\ast$-orbit
with a given property should be interpreted as meaning that
$\ast$-orbits with the given property exist, and all have the same
data.

\begin{exmp}
\label{ex:dataeg}
Let $P\in\cP(2/5,1,A)$ be the periodic orbit depicted in
Fig.~\ref{fig:dataeg}. Then $N_0=N_1=N_3=2$ and $N_2=N_4=1$; and
$A=\{0\}$, $B=\emptyset$, and $C=\{1\}$. The cyclic permutation $\pi$
is given by
\[\stackrel{A}{(0,0)}\to(1,1)\to(3,1)\to
\stackrel{C}{(0,1)}\to(2,0)\to(4,0)\to(1,0)\to(3,0).\]

\begin{figure}[htbp]
\begin{center}
\pichere{0.3}{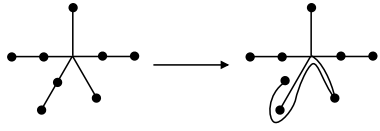}
\end{center}
\caption{An example of a $\ast$-orbit}
\label{fig:dataeg}
\end{figure}
\end{exmp}

The data of a periodic orbit will usually be written in the form
above: the partition $A\cup B\cup C$ is denoted by letters above the
elements $(0,s)$ of $\cL$, and the integers $N_r$ can be deduced from
the elements of $\cL$. For simplicity, it is not explicitly noted that
$(3,0)\to(0,0)$. Note that by $\ast$\,d), the cycle notation of $\pi$
always ends $(m,0)\to(2m,0)\to\cdots\to(n-m,0)$. The reader seeking
further clarification should consult
Examples~\ref{ex:tt},~\ref{ex:simporb}, and~\ref{ex:exorb}, where
other $\ast$-orbits are described in this way.

Write $\tau_1,\tau_2\colon\cL\to\Nset$ for the projections of $\cL$
onto the first and second components respectively; and write
$\pi_1=\tau_1\circ\pi$ and $\pi_2=\tau_2\circ\pi$. Sometimes, with an
abuse of notation, $\pi_1\I$ and $\pi_2\I$ will be used to denote
$\tau_1\circ\pi\I$ and $\tau_2\circ\pi\I$ respectively.

Given $d\in\cD_n$, one can ask whether or not there exists
$P\in\cP(m/n,k,\gamma)$ with $d(P)=d$. This imposes obvious conditions
on $d$, which are expressed by the following result. While hardly
concise, the conditions are in an ideal form for later use.  If
$d\in\cD_n$ satisfies these conditions, it will be said that $d$ is
{\em legal data}.

\begin{lem}
\label{lem:legaldata}
An element $((N_r),\pi,(A,B,C))$ of $\cD_n$ is equal to $d(P)$ for
some $P\in\cP(m/n,k,\gamma)$ if and only if the following conditions
hold:
\begin{enumerate}[\LD  a)]
\item For all $r>0$:
\begin{enumerate}[i)]
\item $\pi_1(r,s)=r+_nm$ for all $s$, and $\pi_2(r,s)$ is
increasing in $s$.
\item $\pi_2(r,0)=0$.
\end{enumerate}
\item For $r=0$:
\begin{enumerate}[i)]
\item $0\in\gamma$.
\item $\pi_1(0,0)=k$.
\item If $c\in\{0,\ldots,N_0-1\}$ is the unique element of $C$, then
$\pi(0,c)=(m,0)$.
\item $\pi_1(0,s)\in\{k,\ldots,m\}$ for all $s$, and is increasing in
$s$.
\item If $s_1<s_2$, $\pi_1(0,s_1)=\pi_1(0,s_2)$ and $s_1\in B\cup C$
then $s_2\in B$.
\item If $s_1<s_2$, $\pi_1(0,s_1)=\pi_1(0,s_2)$, and $s_1,s_2\in A$
(respectively $s_1,s_2\in B$) then $\pi_2(0,s_1)>\pi_2(0,s_2)$
(respectively $\pi_2(0,s_1)<\pi_2(0,s_2)$).
\end{enumerate}
\end{enumerate}
\end{lem}
\smallskip\begin{sketchproof}
The necessity of the conditions is obvious. Their sufficiency can be
shown as follows: divide $e_0$ into $2m+1$ subintervals, each mapped
by $f_{m/n}$ over exactly one edge of $\Gamma_n$; constuct a Markov
graph for $f_{m/n}$ using these subintervals and the edges $e_r$ with
$r>0$ as the Markov partition. Given $k$, $\gamma$, and an element
$d=((N_r),\pi,(A,B,C))$ of $\cD_n$ satisfying \LD a) and \LD b),
conditions a)i), b)iii) and b)iv) ensure that there are exactly two
loops in the Markov graph which are compatible with the first
component $\pi_1$ of $\pi$ and with the partition $(A,B,C)$ (two since
the element $(0,s)$ of $\cL$ with $s\in C$ corresponds to the common
endpoint of two intervals in the partition). Pick either of these
loops, and let $P$ be the periodic orbit of $f_{m/n}$ corresponding to
it. Then conditions b)i) and b)ii) ensure that it belongs to
$\cP(m/n,k,\gamma)$, conditions b)iv)\,--\,vi) ensure that the data of
$P$ agree with the second component $\pi_2$ of $\pi$, and condition
a)ii) ensures that $P$ is a $\ast$-orbit.
\end{sketchproof}\smallskip

\begin{exmp}
\label{ex:simpren}
In this example, Lemma~\ref{lem:legaldata} will be used to show that
there is a bijection
\[\psi\colon\cP(1/3,0,B)\to\cP(1/2,0,B)\]
given by deleting all occurences of $(2,s)$ in the cycle
representation of $\pi$. Thus if
$d(P)= \left((N_0,N_1,N_2),\pi,(A,B,C)\right)$, then
$d(\psi(P))=\left((N_0,N_1),\pi',(A,B,C)\right)$, where $\pi'$ is
obtained from $\pi$ by deleting each occurence of $(2,s)$ in its cycle
representation.

$\psi$ is first shown to be well defined: if
$d(P)=\left((N_0,N_1,N_2),\pi,(A,B,C)\right)$ is legal data, then so
is $\left((N_0,N_1),\pi',(A,B,C)\right)$. Notice that the legality of
$d(P)$ implies that $N_1=N_2$, that $\pi(1,s)=(2,s)$ for all $s<N_1$,
and that $\pi_1(2,s)=0$ for all $s$. Hence
\[\pi'(r,s)=\left\{\begin{array}{ll}
\pi(r,s) & \mbox{ if }r=0\\
\pi^2(r,s) & \mbox{ if }r=1.
\end{array}\right.\]
Now condition~\LD a) holds since $\pi_1'(1,s)=\tau_1\circ\pi^2(1,s)=0$,
$\pi_2'(1,s)=\tau_2\circ\pi^2(1,s)$ is increasing in $s$ by
condition~\LD a)i) for $d(P)$, and
$\pi'(1,0)=\pi^2(1,0)=(0,0)$. Condition~\LD b) holds because $\pi'=\pi$
when $r=0$.

Define $\phi\colon\cP(1/2,0,B)\to\cP(1/3,0,B)$ by following every
$(1,s)$ with $(2,s)$ in the cycle representation of $\pi$. That is, if
$d(P)=\left((N_0,N_1),\pi,(A,B,C)\right)$, then
$d(\phi(P))=\left((N_0,N_1,N_1),\pi',(A,B,C)\right)$, where $\pi'$ is
obtained from $\pi$ by following each occurence of $(1,s)$ with
$(2,s)$ in its cycle notation. It can be shown as above that $\phi$ is
well defined, and it is clearly an inverse to $\psi$.
\end{exmp}
This sort of argument is routine and tedious, and will be abbreviated
when it occurs in earnest in the proof of Theorem~\ref{thm:casem1}.

\subsection{The Train Track Condition}

The main question addressed in this paper is: for which
 \mbox{$P\in\cP(m/n)$} is $\Gamma_n^P$ itself a train track for the
 isotopy class of $F_{m/n}$ relative to $P$? The essential property of
 a train track is {\em efficiency}: intuitively, this says that at any
 point $x$ of $\Gamma_n^P$ at which $(f_{m/n}^P)^i$ is not locally
 injective for some $i$, the image $(f_{m/n}^P)^i(U)$ of a small
 neighbourhood $U$ of $x$ `wraps around' a point of $P$, and hence
 cannot be `pulled tight' without passing through $P$. Efficiency can
 be detected as follows. By condition $\ast$\,d), the only point at
 which such pulling tight could occur is $(m,0)$. For each point $p$
 of $P$ which is not an endpoint of $\Gamma_n^P$ (i.e. $p\not=(r,0)$
 for $0\le r<n$), an arc passing on one side of $p$ will wrap around
 $(m,0)$ under the appropriate iterate of $f_{m/n}^P$, while an arc
 passing on the other side will not, and its image can be pulled
 tight. Provided all of the arcs in the image $f_{m/n}^P(\Gamma_n^P)$
 of $\Gamma_n^P$ pass on the `correct' side of each point of $P$, the
 map $f_{m/n}^P$ is efficient, and hence is a train track map for the
 isotopy class of $F_{m/n}\colon D^2\setminus P\to D^2\setminus P$.

The first part of the next definition gives a partition of
$P\setminus\{(r,0):0\le r<n\}$ into two sets, $\alpha$ and $\beta$:
for those points $p$ of $P$ in $\alpha$, an arc passing $p$ to the
left (according to the orientation of the edge $e_r$ containing $p$)
wraps around $(m,0)$; while for those in $\beta$, an arc passing to
the right wraps around $(m,0)$. It is then possible to make a
combinatorial definition of efficiency, as given by conditions \TT
a)\,--\,\TT d) below.

\begin{defns}
Let $P\in\cP(m/n,k,\gamma)$ have data
$d(P)=\left((N_r),\pi,(A,B,C)\right)$. Define a partition
\[\cL\setminus\{(r,0)\,:\,0\le r<n\}=\alpha\cup\beta\]
inductively as follows. $\pi\I(m,0)\in\beta$. For each $i$ with $2\le
i\le \card P-n$, the two elements $\pi^{-i}(m,0)$ and
$\pi^{-i+1}(m,0)$ are in different sets if $\tau_1(\pi^{-i}(m,0))=0$
and $\tau_2(\pi^{-i}(m,0))\in A$, and are in the same set otherwise.

$P$ is a {\em train track orbit} (or $P$ is {\em TT}) if
and only if the following conditions hold:
\begin{enumerate}[\TT a)]
\item For all $0\le s< N_0-1$,
\[\pi(0,s)\in\left\{\begin{array}{ll}
\alpha&\mbox{ if }s\in B\\
\beta&\mbox{ if }s\in A.
\end{array}\right.\]

\item If $k<r<m$ then $N_{r+n-m}=1$.
\item  If $\pi_2(k+n-m,s)>\pi_2(0,0)$ for some $s>0$,
then $\pi(k+n-m,s)\in\beta$.
\item If $\gamma=A$ then $\pi_2(k+n-m,s)>\pi_2(0,0)$ for all $s>0$.
\end{enumerate}

The set of all train track orbits in $\cP(m/n,k,\gamma)$, in
$\cP(m/n,k)$, and in $\cP(m/n)$ will be denoted $TT(m/n,k,\gamma)$,
$TT(m/n,k)$, and $TT(m/n)$.
\end{defns}

\begin{exmp}
\label{ex:tt}
Let $P\in\cP(1/3,0,B)$ have data
\[
\ov{B}{(0,0)}\to
\ovun{A}{\alpha}{(0,1)}\to
\un{\beta}{(1,1)}\to
\un{\beta}{(2,1)}\to
\ovun{C}{\beta}{(0,2)}\to
(1,0)\to(2,0)\] (see Fig.~\ref{fig:nbtorb}). The partition
$\alpha\cup\beta$ is shown on the cycle representation of $\pi$. This
partition is easy to write down: start with $(0,2)$ (the unique
element of~$C$) in $\beta$, and move backwards through the
permutation, switching from $\beta$ to~$\alpha$ or vice-versa at each
occurence of $A$.  

\begin{figure}[htbp]
\begin{center}
\pichere{0.3}{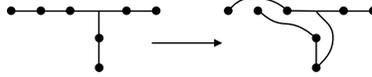}
\end{center}
\caption{A periodic orbit satisfying the TT conditions}
\label{fig:nbtorb}
\end{figure}

The TT conditions can easily be checked:
\begin{enumerate}[a)]
\item $\pi(0,0)\in\alpha$ and $\pi(0,1)\in\beta$.
\item Vacuous, since there are no $r$ with $0=k<r<m=1$.
\item $k+n-m=0+3-1=2$, so the condition requires that
$\pi(2,s)\in\beta$ whenever $\pi_2(2,s)>1$: i.e. $\pi(2,1)\in\beta$,
which is true.
\item Vacuous, since $\gamma=B$.
\end{enumerate}
\end{exmp}

The TT conditions have intuitive motivations. \TT a) ensures that arcs
of $f_{m/n}^P(\Gamma_n^P)$ pass on the `correct' side of images of
points of $P\cap e_0$: for example, if $s\in B$ and $\pi_1(0,s)=r$,
then an arc of $f^P_{m/n}(e_{r+n-m})$ passes to the left of $\pi(0,s)$
and hence $\pi(0,s)$ must lie in $\alpha$ to avoid a violation of
efficiency. \TT b) reflects the fact that for $k<r<m$, the image of
$e_0$ passes images of points of $P\cap e_{r+n-m}$ on both sides:
hence efficiency will be violated if there is any such point other
than the one with image $(r,0)$. \TT c) reflects that the image of
$e_0$ passes to the right of those images of points of $P\cap
e_{k+n-m}$ which are further from the initial point of $e_k$ than
$\pi(0,0)$; and \TT d) is a stronger version of the same condition in
the case $\gamma=A$, when the image of $e_0$ passes on both sides of
points of $e_k$ between the initial point and $\pi(0,0)$.

The TT conditions above are stated in a way which reflects this
intuitive motivation. However they can be replaced by equivalent
alternatives which are easier to work with in practice. This will be
done in Lemma~\ref{lem:equivTT} below.

\begin{figure}[htbp]
\begin{center}
\pichere{0.6}{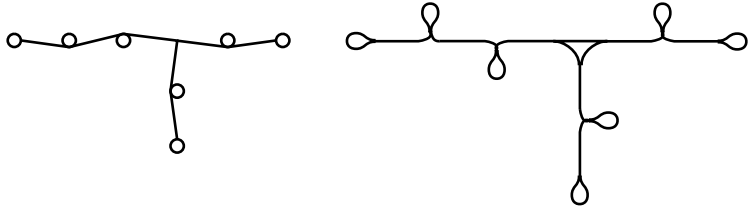}
\end{center}
\caption{Bestvina-Handel and Thurston Train Tracks for the example of Fig.~\ref{fig:nbtorb}}

\label{fig:tracks}
\end{figure}

The reader familiar with this approach to train tracks should accept
without further argument that the TT conditions are the appropriate
combinatorial expression of efficiency. However, a more precise
statement is given by the following definition (which connects the
intuitive notion of a train track with the definition due to Bestvina
and Handel~\cite{BH}, and in particular, by introducing {\em
peripheral loops} around the points of $P$, makes it possible to be
precise about the concept of `wrapping around' a point of $P$) and by
Theorem~\ref{thm:tt}. The rather contorted definition can most easily
be understood by referring to Fig.~\ref{fig:tracks}.

\begin{defns}
Let $P\in\cP(m/n)$, and let $L\cup R$ be a partition of
$\cL\setminus\{(r,0):0\le r<n\}$. The {\em $(L,R)$-Bestvina-Handel
star graph} $BH(P,L,R)\subset T_n$ of $P$ is defined as follows:
\begin{enumerate}[a)]
\item There are $\card P+1$ vertices; one vertex $v$ contained in the 
central decomposition element of $T_n$, and $\card P$ vertices
$\{v_{r,s}:(r,s)\in\cL\}$ with $v_{r,s}$ contained in the same
decomposition element of $T_n$ as the point $(r,s)$ of $P$. If $s>0$
then $v_{r,s}$ is to the left or right of $(r,s)$ (with respect to the
orientation of $e_r$) according as $(r,s)\in L$ or $(r,s)\in R$. (If
$s=0$ then $v_{r,s}$ can be chosen to the left or the right of
$(r,s)$.)
\item There are $\card P$ {\em peripheral} edges
$\{p_{r,s}:(r,s)\in\cL\}$: the peripheral edge $p_{r,s}$ has both
endpoints at $v_{r,s}$, and forms a loop bounding a disk with the
point $(r,s)$ of $P$ in its interior. The peripheral edges are small
enough that no two of them intersect any given decomposition element
of $T_n$.
\item There are $\card P$ {\em main} edges
$\{e_{r,s}:(r,s)\in\cL\}$: the main edge $e_{r,s}$ goes from $v_{r,s}$
to $v_{r,s+1}$ if $s<N_r-1$, and goes from $v_{r,s}$ to $v$ if
$s=N_r-1$. The main edges are chosen so that each interval
decomposition element of $T_n$ contains at most one point of the union
of the interiors of the main edges.
\end{enumerate}
Since $BH(P,L,R)$ is a spine of $T_n\setminus P$, the thick tree map
$F_{m/n}\colon T_n\to T_n$ induces a well defined homotopy class of
graph maps $g\colon BH(P,L,R)\to BH(P,L,R)$. If $g$ is required to
restrict to a homeomorphism of the subgraph of peripheral edges, to
send vertices to vertices, and to be locally injective away from its
vertices, then the image edge-paths of $g$ are also well-defined.
$BH(P,L,R)$ is said to be a {\em Bestvina-Handel star train track} for
$P$ if such a graph map is a train track map: that is, if
\begin{enumerate}[a)]
\item It is {\em absorbed}: the image edge path of each main edge
begins and ends with main edges.
\item It is {\em efficient}: there is no backtracking in the edge path
$g^i(e_{r,s})$ for any $(r,s)\in\cL$ and $i>0$.
\end{enumerate}
The $\ast$-orbit $P\in\cP(m/n)$ is said to {\em have a Bestvina-Handel
star train track} if $BH(P,L,R)$ is a Bestvina-Handel star train track
for some choice of the partition $\cL\setminus\{(r,0):0\le r<n\}=L\cup
R$.
\end{defns}

\begin{rems}
\begin{enumerate}[a)]
\item If $BH(P,L,R)$ is efficient but not absorbed, then there is some
other partition $L'\cup R'$ such that $BH(P,L',R')$ is both efficient
and absorbed (obtained, for example, by applying the Bestvina-Handel
operation of {\em absorbing into the peripheral subgraph}). Requiring
train track maps to be absorbed as well as efficient means that if $P$
has a Bestvina-Handel star train track, then there is a unique choice
of partition $L\cup R$ such that $BH(P,L,R)$ is a Bestvina-Handel star
train track (see the first paragraph of the proof of
Theorem~\ref{thm:tt}). 
\item It follows from the results of~\cite{BH} that if $P$ has a
Bestvina-Handel star train track, and if the transition matrix for the
main edges of that star train track is irreducible, then $P$ has
pseudo-Anosov braid type, and the train track yields a Markov
partition for the pseudo-Anosov representative in the isotopy class of
$F_{m/n}$ in $D^2\setminus P$. Moreover, because the $n$ edge germs at
$v$ are permuted by $g$ they all lie in different gates: hence such a
pseudo-Anosov must have an interior $n$-pronged singularity (whose
prongs are rotated by $m$ under the action of the pseudo-Anosov),
$1$-pronged singularities at each of the points of $P$, and a $(\card
P-n)$-pronged singularity at the boundary.
\end{enumerate}
\end{rems}

\begin{thm}
\label{thm:tt}
Let $P\in\cP(m/n,k,\gamma)$. Then $P$ has a Bestvina-Handel star train
track if and only if $P$ is TT.
\end{thm}

\smallskip\begin{sketchproof}
A straightforward argument shows that $BH(P,L,R)$ is absorbed if and
only if $L=\beta$ and $R=\alpha$: the procedure for constructing
inductively the unique partition $L\cup R$ for which $BH(P,L,R)$ is
absorbed is identical to that by which the partition $\alpha\cup\beta$
is defined.

Thus it is only required to show that $BH(P,\beta,\alpha)$ is
efficient if and only if $P$ is TT. For each $u\not=v$ with $0\le
u,v<n$, write $t_{u,v}$ for the edge-path
$e_{\scriptscriptstyle{u,N_u-1}}\be_{\scriptscriptstyle{v,N_v-1}}$,
and observe that $BH(P,\beta,\alpha)$ is efficient if and only if the
image edge-path $g(e_{r,s})$ of each main edge can be written in the
form
\[g(e_{r,s})=m_1p_1m_2p_2\ldots,m_{k-1}p_{k-1}m_k,\]
where each $p_i$ is a peripheral edge (or its inverse), and each $m_i$
is either a main edge (or its inverse), or one of the edge-paths
$t_{u,v}$. If this condition holds, then a straightforward induction
shows that it holds also for $g^n(e_{r,s})$ for all $n>0$, and hence
there cannot be any cancellation in any of these edge-paths. If the
condition fails, then the edge-path $g(e_{r,s})$ contains some word
$e_{u,v}e_{u,v+1}$ (or its inverse), which under iteration yields the
word $\be_{m,0}e_{m,0}$.

If $P$ is not TT, then one of the conditions \TT a)\,--\,\TT d)
fails. In each case, an argument arising from the intuitive motivation
of the conditions can be used to show that there is some $g(e_{r,s})$
which contains a word of the form $e_{u,v}e_{u,v+1}$ or its
inverse. Hence $BH(P,\beta,\alpha)$ is not efficient.

For the converse, suppose that $P$ is TT. If $BH(P,\beta,\alpha)$ is
not efficient, then there is some edge $e_{u,v}$ such that
$g(e_{u,v})$ contains a word of the form $e_{r,s-1}e_{r,s}$ or
$\be_{r,s}\be_{r,s-1}$. This implies that $u=0$ or $u\ge k+n-m$
(i.e. that $0\le r\le m$), since otherwise
$g(e_{u,v})=e_{u+_nm,v}$ is an edge-path of length $1$.

There are two cases to consider. First suppose that
$\pi_1\I(r,s)\not=0$. Then \TT b) gives $r=k$ (and hence
$\pi_1\I(r,s)=k+n-m$). Thus $u=0$, since only $e_0$ and $e_{k+n-m}$
have images intersecting $e_k$. If $\gamma=B$ then, since $g(e_{0,v})$
passes to the right of $(r,s)$, it can only contain the word
$e_{k,s-1}e_{k,s}$ if $\pi_2(0,v)<s$ and $(r,s)\in\alpha$. This
contradicts \TT c). If $\gamma=A$ then a similar contradiction to \TT
c) and \TT d) arises.

The case where $\pi_1\I(r,s)=0$ can be treated similarly using \TT a).
\end{sketchproof}\smallskip

\begin{exmp}
The Bestvina-Handel and Thurston train tracks for the TT orbit of
Example~\ref{ex:tt} are shown in Fig.~\ref{fig:tracks}. The
Bestvina-Handel train track is obtained by replacing each point of the
orbit with a peripheral loop, and putting the vertex of this loop on
the left or right according as the point of the orbit belongs to
$\beta$ or $\alpha$. The Thurston train track is obtained from the
Bestvina-Handel train track using the techniques of~\cite{BH} (in the
case of $\ast$-orbits, this is simply a matter of replacing the
central valence $n$ vertex with an $n$-gon, and each of the loops with
a $1$-gon, as shown in the figure).
\end{exmp}

As mentioned above, there is an alternative form of the TT conditions
which is often easier to work with:
\begin{lem}
\label{lem:equivTT}
Consider the conditions
\begin{enumerate}[\TTp a)]
\item $\card A=1$, and if $A=\{s\}$ and $j>0$ is least such that
$\tau_1(\pi^j(0,s))=0$, then $C=\{\tau_2(\pi^j(0,s))\}$.
\addtocounter{enumi}{2}
\item If $\gamma=A$ then $\pi(0,0)=(k,1)$.
\end{enumerate}
Then \TT a) and \TTp a) are equivalent; and if they hold then \TT d) and
\TTp d) are equivalent. Moreover, if \TT a) and \TT c) hold then there is
at most one $s$ with $\pi_2(k+n-m,s)\ge\pi_2(0,0)$.
\end{lem}
\begin{proof}
The equivalence of \TT a) and \TTp a) follows easily from the definition
of the partition $\alpha\cup\beta$: note that \TTp a) says simply that
reading the partition $A\cup B\cup C$ along the cycle notation of
$\pi$, one sees $B^{N_0-2}AC$. In particular, since the cycle notation
starts with $(0,0)$, this implies that if $\gamma=A$ then $N_0=2$:
since $\pi_1(0,1)=m$ and $\pi_1(0,0)=k<m$, the equivalence of \TT d) and
\TTp d) follows. For the final statement, note that the elements of
$\beta$ are precisely those which lie between the element of $A$ and
the element and the element of $C$ in the cycle notation of $\pi$, and
only one such can lie in $e_k$: also it is not possible that
$\pi_2(k+n-m,s)=\pi_2(0,0)$, since this would imply
$\pi(k+n-m,s)=\pi(0,0)$.
\end{proof}

If $P\in TT(m/n)$, then it follows from the results of~\cite{BH} that
the set of braid types forced by $\bt(P)$ is precisely the set of braid
types of periodic orbits of $f_{m/n}$ whose span is contained in the
span of $P$. In particular, Lemma~\ref{lem:span} gives:

\begin{thm}
\label{thm:totord}
The set $\{\bt(P,F_{m/n})\,:\, P\in TT(m/n)\}$ is totally ordered
by the forcing relation.
\end{thm}

\subsection{The TT condition in the horseshoe}
The elements of $TT(1/2)=TT(1/2,0,B)$ were calculated in~\cite{Ha}: in
this case the partition $\alpha\cup\beta$ is precisely that defined on
page~880 of~\cite{Ha}, and conditions \TT a) and~\TT c) are equivalent
to the condition that $\pi$ has {\em no bogus transitions}, also on
page~880 of~\cite{Ha} (while \TT b) and~\TT d) are
vacuous). Theorem~2.1 and lemma~2.4 of~\cite{Ha} then give the
following result (in which the words $c_q$ are as defined in
Section~\ref{sec:hsintro}).

\begin{thm}
\label{thm:hstt}
The function $\Qset\cap(0,1/2)\to\cP(1/2,0,B)$ which takes a rational
$q$ to the periodic orbit $P_q$ with code $c_q\opti$ is a bijection
onto $TT(1/2,0,B)$. If $r/s\in(0,1/2)$, the data
$d(P_{r/s})=((N_0,N_1),\pi,(A,B,C))$ of $P_{r/s}$ satisfies
$N_0=s-r+1$ and $N_1=r+1$.
\end{thm}

\subsection{Renomalizing horseshoe TT orbits (the case $k=m-1$)}
\label{sec:casem1}

Theorem~\ref{thm:hstt} gives a relatively straightforward approach to
determining the elements of $TT(m/n,m-1,B)$ for each $m/n\in(0,1/2)$:
the next result defines a `renormalization operator' $\phi$, which is
a bijection from the set $\cP(1/2,0,B)$ of horseshoe periodic orbits
to the set $\cP(m/n,m-1,B)$ with the property that $\phi(P)$ is TT if
and only if $P$ is TT. Example~\ref{ex:simpren} is the simplest
case of this operator, when $m/n=1/3$.

For notational simplicity, $m/n$ will be taken to be a fixed rational
in $(0,1/2)$ throughout this section, and the dependence of some
objects (such as $\phi$) upon it will be dropped.

\begin{thm}
\label{thm:casem1}
There is a bijection $\phi\colon\cP(1/2,0,B)\to\cP(m/n,m-1,B)$ given
by replacing each occurence of $\ov{\gamma}{(0,s)}$ (where
$\gamma\in\{A,B,C\}$) in the cycle representation of $\pi$ by
\[(m-1,s)\to(2m-1,s)\to\cdots\to(n-m,s)\to\ov{\gamma}{(0,s)}\]
and each occurence of $(1,s)$ by
\[(m,s)\to(2m,s)\to\cdots\to(n-m-1,s)\to(n-1,s).\]
Moreover, $\phi$ restricts to a bijection $TT(1/2,0,B)\to TT(m/n,m-1,B)$.

\end{thm}

\begin{exmp}
\label{ex:simporb}
Consider the periodic orbit $P\in\cP(1/2,0,B)$ given by
\[
\ov{B}{(0,0)}\to
\ovun{A}{\alpha}{(0,1)}\to
\un{\beta}{(1,1)}\to
\ovun{C}{\beta}{(0,2)}\to
(1,0)\]
(see Fig.~\ref{fig:simporb}).
\begin{figure}[htbp]
\begin{center}
\pichere{0.6}{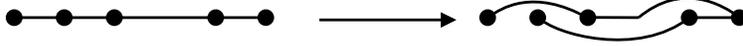}
\end{center}
\caption{A periodic orbit in the horseshoe}
\label{fig:simporb}
\end{figure}

Then $\phi(P)\in\cP(2/5,1,B)$ is given by
\[
\un{\text{from }(0,0)}{\underbrace{\un{\phantom{\alpha}}{(1,0)}\to(3,0)\to
\ov{B}{(0,0)}}}\to
\un{\text{from }(0,1)}{\underbrace{\un{\alpha}{(1,1)}\to
\un{\alpha}{(3,1)}\to
\ovun{A}{\alpha}{(0,1)}}}\to
\un{\text{from }(1,1)}{\underbrace{\un{\beta}{(2,1)}\to
\un{\beta}{(4,1)}}}\]\[\to
\un{\text{from }(0,2)}{\underbrace{\un{\beta}{(1,2)}\to
\un{\beta}{(3,2)}\to
\ovun{C}{\beta}{(0,2)}}}\to
\un{\text{from }(1,0)}{\underbrace{\un{\phantom{\alpha}}{(2,0)}\to(4,0)}}\]
(see Fig.~\ref{fig:renorb}).
It can easily be checked that both of these orbits are TT. The figures
also show the idea of the construction: the pattern of the periodic
orbit of $f_{1/2}$ is replicated in the edges $e_{m-1}$ and $e_m$ of
$\Gamma_n$. 
\begin{figure}[htbp]
\begin{center}
\pichere{0.6}{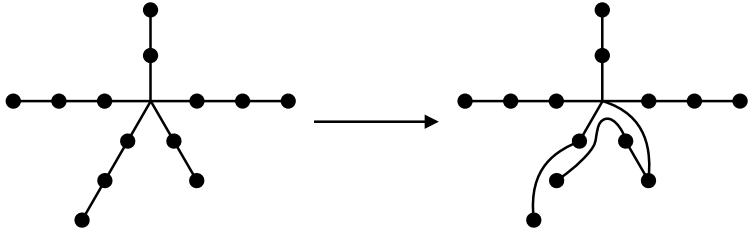}
\end{center}
\caption{The periodic orbit in $\cP(2/5,1,B)$ obtained from that of Fig.~\ref{fig:simporb}}
\label{fig:renorb}
\end{figure}
\end{exmp}

\begin{proof}
The first step is to show that the construction yields a well-defined
function $\phi\colon\cP(1/2,0,B)\to\cP(m/n,m-1,B)$: that is, that it
transforms legal data to legal data.

Suppose, then, that $P\in\cP(1/2,0,B)$ has data
$d(P)=((N_0,N_1),\pi,(A,B,C))$. Since $\cO_{m/n}[m-1,0]$ and
$\cO_{m/n}[m,n-1]$ partition $\Zset_n$, if the construction does give
an element $P'$ of $\cP(m/n,m-1,B)$ it must have data
\[d(P')=\left((N'_r),\pi',(A,B,C)\right),\]
where
\[N'_r=\left\{\begin{array}{ll}
N_0&\text{ if }r\in\cO_{m/n}[m-1,0]\\
N_1&\text{ if }r\in\cO_{m/n}[m,n-1],
\end{array}\right.\]
and
\[\pi'(r,s)=\begin{cases}
(\pi_1(0,s)+m-1,\pi_2(0,s)) & \text{ if }r=0\\
(\pi_1(1,s)+m-1,\pi_2(1,s)) & \text{ if }r=n-1\\
(m+_nr,s) & \text{ otherwise.}
\end{cases}\]
(The sets $A$, $B$, and $C$ are unchanged by construction.)

It is routine to show that conditions \LD a) and \LD b) of
Lemma~\ref{lem:legaldata} hold for this data, and hence that
$\phi\colon\cP(1/2,0,B)\to\cP(m/n,m-1,B)$ is well defined as
required. To show that it is a bijection, observe that if
$P\in\cP(m/n,m-1,B)$ has data $d(P)=((N_r),\pi,(A,B,C))$, then
conditions \LD a)i) and \LD b)iv) ensure that every occurrence of
$(m-1,s)$ in the cycle representation of $\pi$ is followed by
\[(m-1,s)\to(2m-1,s)\to\cdots\to(n-m,s)\to\ov{\gamma}{(0,s)}\]
for some $\gamma\in\{A,B,C\}$, and every occurrence of $(m,s)$ is
followed by
\[(m,s)\to(2m,s)\to\cdots\to(n-m-1,s)\to(n-1,s).\]
Replacing the first type of block with $\ov{\gamma}{(0,s)}$ and the
second with $(1,s)$ thus defines an inverse
$\psi\colon\cP(m/n,m-1,B)\to\cP(1/2,0,B)$ of $\phi$ (checking that the
data so obtained is legal is again routine).

Finally, it must be shown that the bijection $\phi$ preserves the TT
conditions.  Let $P\in\cP(1/2,0,B)$ and $P'=\phi(P)\in\cP(m/n,m-1,B)$.
Since the sets $A$, $B$, and $C$ are the same for $P$ and $P'$, the
condition \TTp a) is satisfied either for neither or for both of $P$
and $P'$. Conditions \TT b) and \TT d) are vacuous for both $P$ and
$P'$. Condition \TT c) for $P$ reads
\[\pi(1,s)\in\beta\text{ whenever }\pi_2(1,s)>\pi_2(0,0),\]
while for $P'$ it reads
\[\pi'(n-1,s)\in\beta\text{ whenever }\pi'_2(n-1,s)>\pi'_2(0,0).\]
Since $\pi'(n-1,s)\in\beta$ if and only if $\pi(1,s)\in\beta$;
$\pi'_2(n-1,s)=\pi_2(1,s)$; and $\pi'_2(0,0)=\pi_2(0,0)$ (these last
two by the expression for $\pi'$ above), \TT c) is also satisfied
either for neither or for both of $P$ and $P'$.
\end{proof}

The next step is to use Lemma~\ref{lem:identify} to identify the
orbits $\phi(P_{r/s})$ which make up $TT(m/n,m-1,B)$. The function
$\xi_{u/v,p/q}$ in the statement of the next theorem is the one
defined in Lemma~\ref{lem:xi}.

\begin{thm}
\label{thm:identifym1}
Let $u/v=LFP(m/n)$ and $p/q=RFP(m/n)$.  Let $r/s\in(0,1/2)$. Then
$\phi(P_{r/s})$ has height $\xi_{u/v,p/q}(r/s)$ and decoration
$w_{m/n}$.
\end{thm}
\begin{proof}
By Theorem~\ref{thm:hstt}, $\phi(P_{r/s})$ has period
\[(s-r+1)\card\cO_{m/n}[m-1,0]+(r+1)\card\cO_{m/n}[m,n-1],\]
which, by Lemma~\ref{lem:fareyprops}, is equal to
$(s-r+1)v+(r+1)q$. It will be shown that $\phi(P_{r/s})$ has rotation
interval
\[\left[\frac{(s-r)u+rp}{(s-r)v+rq},\frac{u+p}{v+q}\right]
=\left[\xi_{u/v,p/q}\left(\frac{r}{s}\right),\frac{m}{n}\right],\]
which will establish the result by Lemmas~\ref{lem:identify}
and~\ref{lem:starhs}.

If $P$ is a $\ast$-orbit with data $d(P)=((N_r),\pi,(A,B,C))$, then
the rotation number of $P$ about the fixed point is given by the
number of times it cycles around $\Gamma_n$ divided by its period:
that is,
\[\rho(P)=\frac{\sum_{r=n-m}^{n-1}N_r}{\sum_{r=0}^{n-1}N_r}.\] 

Let $P\in\cP(1/2,0,B)$. Then $\rho(P)=N_1/(N_0+N_1)$,
and if $\phi(P)$ has data $((N'_r),\pi,(A,B,C))$ then by
Lemma~\ref{lem:fareyprops}, and the fact that $N'_r=N_0$ for all
$r\in\cO_{m/n}[m-1,0]$ and $N'_r=N_1$ for all $r\in\cO_{m/n}[m,n-1]$,
\[\frac{\displaystyle{\sum_{r\in\cO_{m/n}[m,n-1],r\ge
n-m}N'_r}}{\displaystyle{\sum_{r\in\cO_{m/n}[m,n-1]}N'_r}}
=\frac{p}{q}\quad\text{ and }\quad
\frac{\displaystyle{\sum_{r\in\cO_{m/n}[m-1,0],r\ge
n-m}N'_r}}{\displaystyle{\sum_{r\in\cO_{m/n}[m-1,0]}}N'_r}=\frac{u}{v},\]
and hence
\[\rho(\phi(P))=\frac{N_0u+N_1p}{N_0v+N_1q}.\]
If $P$ has period $b$ and rotation number $a/b$ then $N_0=b-a$ and
$N_1=a$, and so
\[\rho(\phi(P))=\frac{(b-a)u+ap}{(b-a)v+aq}=\xi_{u/v,p/q}(a/b).\]

Now the rotation interval $[r/s,1/2]$ of $P_{r/s}$ is the set of
rotation numbers of periodic orbits of $f_{1/2}$ contained in the span
of $P_{r/s}$, and, by a theorem of Boyland~\cite{Boyri}, for each such
rotation number $a/b$ there exists such an orbit with period
$b$. Likewise, the rotation interval of $\phi(P_{r/s})$ is the set of
rotation numbers of periodic orbits of $f_{m/n}$ contained in the span
of $\phi(P_{r/s})$; but these orbits are precisely the images under
$\phi$ of those defining the rotation interval of $P_{r/s}$. Hence the
rotation interval of $\phi(P_{r/s})$ is
\[\left\{\xi_{u/v,p/q}\left(\frac{a}{b}\right)\colon 
\frac{a}{b}\in\left[\frac{r}{s},\frac{1}{2}\right]\right\}
=\xi_{u/v,p/q}\left(\left[\frac{r}{s},\frac{1}{2}\right]\right)
=\left[\xi_{u/v,p/q}\left(\frac{r}{s}\right),\frac{m}{n}\right]\]
as required.
\end{proof}

\begin{rem}
If $n\ge 2$ then $\cP(1/n)=\cP(1/n,0,B)$, and hence the results of
this section give a complete description of $TT(1/n)$: there is a
bijection $q\mapsto P_q^{1/n}$ from $(0,1/n)\cap\Qset$ to $TT(1/n)$
such that $P_q^{1/n}$ has height $q$ and decoration
$w_{1/n}=0^{n-3}$. In other words, $P_q^{1/n}$ has the same braid type
as the horseshoe orbits of code $c_q\opti0^{n-3}\opti$.

In general, since $\xi_{u/v,p/q}:(0,1/2)\to(u/v,m/n)$ is an increasing
bijection, the results of this section yield a bijection $q\mapsto
P_q^{m/n}$ from $(\LFP(m/n),m/n)\cap\Qset$ to $TT(m/n,m-1,B)$ with the
property that $P_q^{m/n}$ has height $q$ and decoration $w_{m/n}$.
\end{rem}

\subsection{Admissible $k$}
\label{sec:adm}
Condition \TT b) gives restrictions on the values of $k$ for which
$TT(m/n,k)$ is non-empty: in this section it will be shown that the
number of such {\em admissible} values of $k$ is equal to the length
of the left Farey sequence of $m/n$.

\begin{defn}
Let $m/n\in(0,1/2)$. An integer $k\in[0,m-1]$ is {\em
$m/n$-admissible} if
\[[k+1,m-1]\cap\cO_{m/n}[k,0]=\emptyset.\] The set of all $m/n$-admissible
integers $k$ is denoted $\cA_{m/n}$.
\end{defn}

\begin{lem}
\label{lem:admiss}
Let $m/n\in(0,1/2)$ and $k\in[0,m-1]$. If $k\not\in\cA_{m/n}$, then
$TT(m/n,k)$ is empty.
\end{lem}
\begin{proof}
Suppose $P\in TT(m/n,k)$ has data $d(P)=((N_r),\pi,(A,B,C))$. Then by
\TT b), $N_{r+n-m}=1$ for all $r$ with $k<r<m$. Applying \LD a)\,i)
inductively, it follows that $N_s=1$ for all
$s\in\bigcup_{k<r<m}\cO_{m/n}[m,r+n-m]$. However $N_k>1$, since
$\pi_2(0,0)=\pi_2(k+n-m,0)=k$. Thus $TT(m/n,k)$ must be empty unless
$k\not\in\bigcup_{k<r<m}\cO_{m/n}[m,r+n-m]$. However, this condition
is equivalent to $k\in\cO_{m/n}[r,0]$ for all $r$ with $k<r<m$; which
is in turn equivalent to $r\not\in\cO_{m/n}[k,0]$ for all $r$ with
$k<r<m$: that is, to $m/n$-admissibility.
\end{proof}

The remainder of this section is devoted to studying the structure of
the set $\cA_{m/n}$. It will be seen in Section~\ref{sec:m1} that
the converse to Lemma~\ref{lem:admiss} is also true: if
$k\in\cA_{m/n}$, then $TT(m/n,k)$ is non-empty.

Notice that $0,m-1\in\cA_{m/n}$ for all $m/n$: in particular
$\cA_{1/n}=\{0\}$. Lemma~\ref{lem:constructadm} below connects
$\cA_{m/n}$ with $\cA_{\LFP(m/n)}$, and hence yields an
inductive description of $\cA_{m/n}$ for all $m/n$.

Let $m/n\in(0,1/2)$ with $\LFP(m/n)=u/v$, and write
$R_{m/n}=\cO_{m/n}[m,n-1]$. By Lemma~\ref{lem:fareyprops}, $R_{m/n}$
has cardinality $n-v$. Moreover, $R_{m/n}\cap\cA_{m/n}=\emptyset$: for
if $k\in R_{m/n}$, then $\cO_{m/n}[k,0]$ contains the complement of
$R$, and in particular contains $m-1$.  Define a function
$\psi_{m/n}\colon\Zset_n\to\Zset_v$ by
\[\psi_{m/n}(k)=k-\card \left(R_{m/n}\cap[0,k]\right).\]

Notice that
\[\psi_{m/n}(k)=\left\{
\begin{array}{ll}
\psi_{m/n}(k-1)+1 & \text{ if }k\not\in R_{m/n}\\
\psi_{m/n}(k-1)   & \text{ if }k\in R_{m/n}
\end{array}
\right.\]
It can easily be seen that $\psi_{m/n}(0)=0$,
$\psi_{m/n}(m-1)=\psi_{m/n}(m)=u$ (using Lemma~\ref{lem:fareyprops}),
and $\psi_{m/n}(n-1)=v-1$. In particular, $\psi_{m/n}$ is an
increasing surjection.

\begin{lem}
\label{lem:constructadm}
Let $m/n\in(0,1/2)$ with $m>1$, $u/v=\LFP(m/n)$, and
$k\in[0,m-2]\setminus R_{m/n}$. Then $k\in\cA_{m/n}$ if and only if
$\psi_{m/n}(k)\in\cA_{u/v}$.
\end{lem}
\begin{proof}
Observe that $\psi_{m/n}\vert_{\Zset_n\setminus
R_{m/n}}\colon\Zset_n\setminus R_{m/n}\to\Zset_v$ is a bijection which
conjugates the function $+_{u/v}\colon\Zset_v\to\Zset_v$ given by
$+_{u/v}(k)=k+_vu$ to the function $\tilde+_{m/n}\colon\Zset_n\setminus
R_{m/n}\to\Zset_n\setminus R_{m/n}$ given by
\[\tilde+_{m/n}(k)=\left\{
\begin{array}{ll}
k+_nm&\text{ if }k\not=0\\
m-1&\text{ if }k=0.
\end{array}
\right.
\]
Since $k\in[0,m-2]$ and $m-1\not\in R_{m/n}$, it follows that
$\psi_{m/n}(k)\le\psi_{m/n}(m-1)-1=u-1$. Now
\begin{eqnarray*}
k\in\cA_{m/n}&\iff&r\not\in\cO_{m/n}[k,0]\text{ for }k<r<m\\
&\iff&r\not\in\cO_{\tilde+_{m/n}}[k,0]\text{ for }k<r<m\text{ and
}r\not\in R_{m/n}\\
&\iff&\psi_{m/n}(r)\not\in\cO_{u/v}[\psi_{m/n}(k),\psi_{m/n}(0)]\text{
for }k<r<m\text{ and }r\not\in R_{m/n}\\
&\iff&s\not\in\cO_{u/v}[\psi_{m/n}(k),0]\text{ for
}\psi_{m/n}(k)<s<\psi_{m/n}(m)=u\\
&\iff&\psi_{m/n}(k)\in\cA_{u/v}
\end{eqnarray*}
(where the penultimate equivalence uses that $R_{m/n}$
contains neither $k$ nor $m-1$).
\end{proof}

It follows that $\cA_{m/n}=\{m-1\}\cup\psi_{m/n}\vert_{\Zset_n\setminus
R_{m/n}}\I(\cA_{u/v})$, and in particular that $\cA_{m/n}$ has one
more element than $\cA_{u/v}$. Since $\cA_{1/n}=\{0\}$ for all $n$,
the following result holds:
\begin{cor}
\label{cor:adm}
Let $m/n\in(0,1/2)$. Then $\cA_{m/n}$ has the same cardinality as
$\LFS(m/n)$.
\end{cor}
The following result will be important in the next section.
\begin{lem}
\label{lem:Tcard}
Let $\LFS(m/n)=(0=u_1/v_1,u_2/v_2,\ldots,u_\alpha/v_\alpha)$, and let
the elements of $\cA_{m/n}$ be $0=k_1<k_2<\cdots<k_\alpha=m-1$. Then
$\cO_{m/n}[k_i,0]$ has cardinality $v_i$ for each~$i$. In particular,
the elements of $\cA_{m/n}$ appear in decreasing order along the orbit
$\cO_{m/n}[m,0]$.
\end{lem}
\begin{proof}
That $\cO_{m/n}[m-1,0]$ has cardinality $v_\alpha$ follows from
Lemma~\ref{lem:fareyprops}. The proof of Lemma~\ref{lem:constructadm}
shows that if $k\in\cA_{m/n}$ is distinct from $m-1$, then
$\cO_{u/v}[\psi_{m/n}(k),0]=\psi_{m/n}(\cO_{m/n}[k,0])$, and the
result follows by induction.
\end{proof}

\subsection{Train track orbits with $k<m-1$}
\label{sec:m1}
Let $m/n\in(0,1/2)$ with $m>1$, and
$\LFS(m/n)=(0=u_1/v_1,u_2/v_2,\ldots,u_\alpha/v_\alpha)$. By
Corollary~\ref{cor:adm}, $\cA_{m/n}$ has $\alpha$ elements, which will
be denoted $0=k_1<k_2<\cdots<k_\alpha=m-1$. Let $i<\alpha$, and define
\begin{eqnarray*}
R&=&\cO_{m/n}[m,k_{i+1}+n-m]\\
S&=&\cO_{m/n}[k_{i+1},k_i+n-m]\\
T&=&\cO_{m/n}[k_i,0].
\end{eqnarray*}
The dependence of the sets on $m/n$ and $i$ is suppressed, since these
two variables will remain fixed throughout the section. $R$, $S$, and
$T$ are mutually disjoint by Lemma~\ref{lem:Tcard}, and hence define a
partition of $\Zset_n$. The cardinalities of $T$ and $S$ are given by
Lemma~\ref{lem:Tcard}, and the cardinality of $R$ can therefore be
deduced:
\begin{eqnarray*}
\card R &=& n-v_{i+1}\\
\card S &=& v_{i+1}-v_i\\
\card T &=& v_i.
\end{eqnarray*}

The following theorem gives a complete description of $TT(m/n,k_i,B)$:
\begin{thm}
\label{thm:findperm}
There is a bijection from $(0,1)\cap\Qset$ to $TT(m/n,k_i,B)$, defined
as follows: if $p/q\in(0,1)$, then the element $P_{i,p/q,m/n}$ of
$TT(m/n,k_i,B)$ corresponding to $p/q$ has data
$d(P_{i,p/q,m/n})=((N_r),\pi,(A,B,C))$ given by: $A=\{p\}$,
$B=\{0,\ldots,q-1\}\setminus\{p\}$, $C=\{q\}$,
\[
N_r=\left\{
\begin{array}{ll}
1 & \text{ if }r\in R\\
q+1-p & \text{ if }r\in S\\
q+1 & \text{ if }r\in T,
\end{array}
\right.
\]
and 
\begin{eqnarray*}
\pi(0,s)&=&\left\{
\begin{array}{ll}
(k_i,q-p+s)&\text{ if }0\le s<p\\
(k_{i+1},q-p)&\text{ if }s=p\\
(k_{i+1},s-p)&\text{ if }p+1\le s<q\\
(m,0)&\text{ if }s=q
\end{array}
\right.\\
\pi(k_i+n-m,q-p)&=&(k_i,q)\\
\pi(r,s)&=&(r+_nm,s)\quad\text{ for all other }(r,s).
\end{eqnarray*}
\end{thm}

\begin{figure}[htbp]
\psfrag{1}[r]{$e_{k_i}$}
\psfrag{2}[tl]{$e_{k_{i+1}}$}
\psfrag{3}[t]{$e_m$}
\begin{center}
\pichere{0.4}{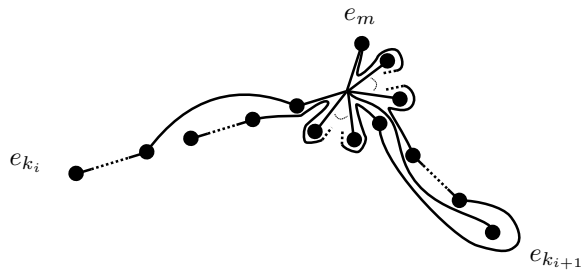}
\end{center}
\caption{A schematic representation of an element of $TT(m/n,k_i,B)$}
\label{fig:scheme}
\end{figure}

A schematic representation of $P_{i,p/q,m/n}$ is shown in
Fig.~\ref{fig:scheme}, which depicts the image of $e_0$ and of
$e_{k_i+n-m}$ through to $e_{n-1}$. $e_0$ contains $q+1$ points of
$P$, of which $p$ are mapped to $e_{k_i}$ (the other $q+1-p$ points of
$P$ on $e_{k_i}$ are images of points of $P$ on $e_{k_i+n-m}$); $q-p$
are mapped to $e_{k_{i+1}}$ in the configuration shown; and one is
mapped to $e_m$.

\begin{figure}[htbp]
\begin{center}
\pichere{0.6}{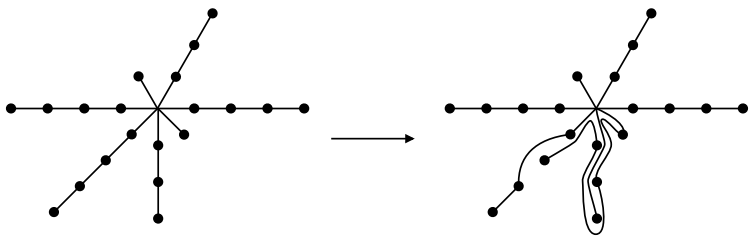}
\end{center}
\caption{The periodic orbit $P_{1,1/3,3/7}$}
\label{fig:exorb}
\end{figure}

\begin{exmp}
\label{ex:exorb}
Let $m/n=3/7$. Now $\LFS(3/7)=(0,1/3,2/5)$ has three elements, and
hence $\cA_{3/7}=\{0,1,2\}$, with $k_1=0$, $k_2=1$, and $k_3=2$. Pick
$i=2$. Then $R=\cO_{3/7}[3,2+7-3]=\{3,6\}$,
$S=\cO_{3/7}[2,1+7-3]=\{2,5\}$, and $T=\cO_{3/7}[1,0]=\{1,4,0\}$. 

Consider the orbit $P_{1,1/3,3/7}$: by the statement of the theorem
(with $p/q=1/3$) it has data
\[
\ov{B}{(0,0)}\to(1,2)\to(4,2)\to
\ov{B}{(0,2)}\to(2,1)\to(5,1)\to(1,1)\]\[\to(4,1)\to
\ov{A}{(0,1)}\to(2,2)\to(5,2)\to(1,3)\to(4,3)\to
\ov{C}{(0,3)}\]\[\to\ov{\phantom{C}}{(3,0)}\to(6,0)\to(2,0)\to(5,0)\to(1,0)\to(4,0)
\]
(see Fig.~\ref{fig:exorb}).
\end{exmp}

\begin{proof}
Suppose $P\in TT(m/n,k_i,B)$ has data $d(P)=((N_r),\pi,(A,B,C))$. Let
$\mu\le N_0-1$ be greatest such that either $\mu\in A$ or
$\pi_1(0,\mu)<k_{i+1}$. Then $\mu\in[1,N_0-2]$ since $0\in B$ and
$\pi_1(N_0-1,0)=m$. It will be shown that the data of $P$ must be
as given in the statement of the theorem, with $p/q=\mu/(N_0-1)$. The
proof is broken down into several short steps.
\begin{enumerate}[i)]
\item $T\cap[k_i+1,m-1]=\emptyset$ (by the definition of
$m/n$-admissibility of $k_i$).
\item If $r\in R$, then $N_r=1$ and $\pi(r,0)=(r+_nm,0)$.

Since $k_i<k_{i+1}<m$, \TT b) gives $N_{k_{i+1}+n-m}=1$, and applying
\LD a)\,i) inductively gives the result.
\item If $r\in T$ then $N_r=N_0$, and if $r\not=0$ then
$\pi(r,s)=(r+_nm,s)$ for all $s<N_0$.

By~i), if $j\in T$ with $j\not=0$, the points of $P\cap e_{j+_nm}$ are
precisely the images of the points of $P\cap e_j$.
\item $\pi(k_i+n-m,N_{k_i+n-m}-1)=(k_i,N_0-1)$.

$N_m=1$ by~ii), so $\pi(0,N_0-1)=(m,0)$. Hence $N_0-1\in\beta$ (by
definition of the partition $\alpha\cup\beta$), and it follows
from~iii) that $(r,N_0-1)\in\beta$ for all $r\in T$: in particular,
$(k_i,N_0-1)\in\beta$. By \TT a), if $\pi_1\I(k_i,N_0-1)=0$ then
$\pi_2\I(k_i,N_0-1)\in A$, contradicting \LD b)\,i),ii),v). Hence
$\pi_1\I(k_i,N_0-1)=k_i+n-m$, and the result follows.
\item If $s<N_{k_i+n-m}-1$, then $\pi(k_i+n-m,s)=(k_i,s)$ (by~iv) and
Lemma~\ref{lem:equivTT}).
\item $S\cap[k_i,m]=\{k_{i+1}\}$.

Suppose not: let $k$ be the last element in the orbit segment from
$k_{i+1}$ to $k_i+n-m$ under addition of $m$ modulo $n$ which lies in
$[k_i,m]$. Then $\cO_{m/n}[k,0]$ contains no elements of $[k+1,m-1]$,
and hence $k\in\cA_{m/n}$; moreover, $k_i<k<k_{i+1}$ by
Lemma~\ref{lem:Tcard}.  This is a contradiction, since $k_i$ and
$k_{i+1}$ are successive elements of $\cA_{m/n}$.
\item $[k_i+1,m]\setminus\{k_{i+1}\}\subseteq R$ (by~i) and~vi)).
\item $A=\{\mu\}$ and $\pi_1(0,s)$ is equal to $k_i$ if $0\le s<\mu$,
to $k_{i+1}$ if $\mu\le s<N_0-1$, and to $m$ if $s=N_0-1$.

That $\pi_1(0,s)$ is equal to $k_i$, $k_{i+1}$ or $m$ for all $s$
follows from~ii) and~vii). Then $A=\{\mu\}$ by the definition of
$\mu$,~ii), and \LD b)\,i),iii),v): the stated values of $\pi_1(0,s)$
are then immediate.
\item If $r\in S$ then $N_r=N_0-\mu$, and if $r\not=k_i+n-m$ then
$\pi(r,s)=(r+_nm,s)$ for all $s<N_0-\mu$.

By~viii), $N_{k_{i+1}}=1+\card[\mu,N_0-2]=N_0-\mu$. The result follows
by~vi) and \LD a)\,i).
\item $\pi(0,\mu)=(k_{i+1},N_0-\mu-1)$

By ii), iii), iv) and ix), there is a segment of the orbit of $\pi$ as
follows: 
\begin{eqnarray*}
(0,s)&\to&\un{S}{\underbrace{(k_{i+1},N_0-\mu-1)\to\cdots\to(k_i+n-m,N_0-\mu-1)}
}\qquad\qquad\\
&\to&\un{T}{
\underbrace{(k_i,N_0-1)\to\cdots\to\ov{C}{(0,N_0-1)}}}\to(m,0)
\end{eqnarray*}
for some $s$. However $s\in A$ by \TTp a), and hence $s=\mu$ by viii).
\end{enumerate}

These results are enough to show that $d(P)$ is as given in the
statement of the theorem with $p=\mu$ and $q=N_0-1$. Conversely, the
data given in the statement defines an element of $TT(m/n,k_i,B)$ by
construction, provided only that $\pi$ is a cyclic permutation. The
proof is therefore completed by showing that given two integers $p$
and $q$ with $0<p<q$, the permutation $\pi$ given in the statement is
cyclic if and only if $p$ and $q$ are coprime. It is clear that for
any $(r,s)$, there is some $K$ such that $\tau_1(\pi^K(r,s))=0$; and
that $\pi^{v_i}(0,0)=(0,q-p)\not=(0,0)$. Hence the cyclicity of $\pi$
is equivalent to the cyclicity of the first return permutation of
$\pi$ on $\{(0,s)\colon 1\le s\le q\}$. Let $\rho$ be this first
return permutation, i.e. $\rho(s)=\tau_2(\pi^{K_s}(0,s))$, where
$K_s>0$ is least such that $\tau_1(\pi^{K_s}(0,s))=0$ and
$\tau_2(\pi^{K_s}(0,s))\not=0$. Then a straightforward calculation
shows that $\rho$ is a rotation of $[1,q]$ by $-p$, and hence is
cyclic if and only if $p$ and $q$ are coprime.
\end{proof}

The next step is to identify the horseshoe braid types of these
orbits: it will be shown that $P_{i,p/q,m/n}$ has period
$n+v_ip+v_{i+1}(q-p)$ and rotation interval
$[\frac{u_ip+u_{i+1}(q-p)}{v_ip+v_{i+1}(q-p)},\frac{m}{n}]=
[\xi_{u_i/v_i,u_{i+1}/v_{i+1}}(p/q),m/n]$, and hence has height
$\xi_{u_i/v_i,u_{i+1}/v_{i+1}}(p/q)$ and decoration $w_{m/n}$ by
Lemma~\ref{lem:identify}. In constrast to the situation in
Section~\ref{sec:casem1}, there is no renormalization operator to
provide a short cut to these rotation intervals: instead, they will be
calculated by Markov partition techniques. Although the calculation is
rather complicated the techniques are quite standard, and as such only
a sketch proof, outlining the main steps, is given, so as to enable
the enthusiastic reader to reconstruct the proof without too much
difficulty.

\begin{thm}
\label{thm:identify}
$P_{i,p/q,m/n}\in TT(m/n,k_i,B)$ has height
$\xi_{u_i/v_i,u_{i+1}/v_{i+1}}(p/q)$ and decoration $w_{m/n}$.
\end{thm}
\smallskip\begin{sketchproof}
That $P=P_{i,p/q,m/n}$ has period $n+v_ip+v_{i+1}(q-p)$ is immediate
from Theorem~\ref{thm:findperm} and the cardinalities of the sets $R$,
$S$, and $T$. The rotation interval of $P$ can be determined using
Markov partition techniques, using the partition of $\Gamma_n^P$ into
intervals whose endpoints are the points of $P$ and the valence $n$
vertex~$v$: the interval with endpoints $(r,s)$ and either $(r,s+1)$
or $v$ is labelled $<r,s>$. The Markov graph with these intervals as
vertices can in principle be determined from the expression for $\pi$
given in Theorem~\ref{thm:findperm}. Each loop in the Markov graph
corresponds to a periodic orbit $Q$ of $f_{m/n}$ whose braid type is
forced by that of $P$: the rotation interval of $P$ is the set of
rotation numbers of such orbits about $v$. To calculate the
rotation number corresponding to a given loop, one counts the number
of times it goes around the star (i.e. the number of occurences of
$<r,s>$ in the loop where $r\ge n-m$), and divides by its
length. 

It is clear that no such orbit can have rotation number greater than
$m/n$, and that there is an orbit with this rotation number (namely
the one given by the loop through the intervals with endpoint $v$).
Thus it only remains to calculate the smallest possible rotation
number; by standard arguments, this will be realized by a {\em
minimal} loop (i.e. one which passes through each vertex at most
once). The full Markov graph is too complicated to study in its
entirety, so a sequence of simplifications is made.
\begin{enumerate}[i)]
\item Every loop must contain $<0,s>$ for some $s$ (i.e. $\{<0,s>:0\le
s\le q\}$ is a {\em rome}~\cite{BGMY} for the Markov graph): the
structure of all loops is therefore given by considering only these
intervals, and making a list of minimal paths between such
intervals. Each such {\em basic} path has associated a rotation number
(the number of occurences of $<r,s>$ with $r\ge n-m$ divided by its
length), and the rotation number of a loop made by concatenating basic
paths is the Farey sum of the corresponding rotation numbers.

Because of the structure of $\pi$ (the only transition from
$<r,s>$ is to $<r+_nm,s>$ unless $r=0$, $k_i+n-m$, or $k_{i+1}+n-m$),
most of these paths correspond to passing either through $T$, or
through $S\cup T$, or through $R\cup S\cup T$: a short calculation
(similar in spirit to the proof of Lemma~\ref{lem:fareyprops}) shows
that the rotation numbers corresponding to these three types of path
are $u_i/v_i<u_{i+1}/v_{i+1}<m/n$ respectively.
\item There are basic paths from $<0,q>$ to $<0,s>$ for all $s\in[0,q]$,
each with rotation number $m/n$. Any loop passing through $<0,q>$ can
be replaced by one with smaller rotation number, and $<0,q>$ can
therefore be ignored.
\item There is a {\em rotation} loop which passes through 
each $<0,s>$ with $0\le s<q$ exactly once, made of basic paths from
$<0,s>$ to $<0,s-_qp>$, each of which has rotation number either
$u_i/v_i$ or $u_{i+1}/v_{i+1}$. Hence given any loop which uses a
basic path with rotation number greater than $u_{i+1}/v_{i+1}$,
another loop of smaller rotation number can be constructed by
replacing this basic path with a segment of the rotation loop. Thus
the loop with minimal rotation number cannot use basic paths with
rotation numbers greater than $u_{i+1}/v_{i+1}$.
\item The only basic paths which now remain to be considered are the
following:
\[
<0,s>\to<0,q-p+s>\text{ for $0\le s\le p-1$, with rotation
$u_i/v_i$} \]
\[<0,p>\to<0,t>\text{ for $0\le t\le q-1$, with rotation
$u_{i+1}/v_{i+1}$}\]
\[<0,s>\to<0,s-p>\text{ for $p+1\le s\le q-2$, with rotation
$u_{i+1}/v_{i+1}$}\]
\[<0,q-1>\to<0,t>\text{ for $q-p-1\le t\le q-1$, with rotation
$u_{i+1}/v_{i+1}$.} \]

One can then argue that basic paths from $<0,p>$ to $<0,t>$ for
$t\not=0$, or from $<0,q-1>$ to $<0,t>$ for $t\not=q-p-1$ could be
replaced by compound paths with lower rotation number. Hence the
rotation loop realizes the minimum rotation number. 
\end{enumerate}
In the rotation loop, $p$ of the basic paths have rotation number
$u_i/v_i$, and the other $q-p$ have rotation number
$u_{i+1}/v_{i+1}$. Hence $P_{i,p/q,m/n}$ has rotation interval
\[\left[\frac{u_ip+u_{i+1}(q-p)}{v_ip+v_{i+1}(q-p)},\frac{m}{n}\right]
=\left[\xi_{u_i/v_i,u_{i+1}/v_{i+1}}(p/q),\frac{m}{n}\right].\]
Since it has period $n+v_ip+v_{i+1}(q-p)$, the result follows by
Lemma~\ref{lem:identify}.
\end{sketchproof}\smallskip

Hence $TT(m/n,k_i,B)$ consists of exactly one orbit with height $r/s$
and decoration $w_{m/n}$ for each
$r/s\in(u_i/v_i,u_{i+1}/v_{i+1})$. Combining this with the results of
Section~\ref{sec:m1}, $\bigcup_kTT(m/n,k,B)$ consists of one orbit
with height $r/s$ and decoration $w_{m/n}$ for each $r/s\in
(0,m/n)\setminus\LFS(m/n)$. The `missing' heights
$r/s\in\LFS(m/n)\setminus\{0\}$ are supplied by orbits in
$\bigcup_kTT(m/n,k,A)$, as described by Theorem~\ref{thm:caseA}
below. The case $\gamma=A$ is much easier than the case $\gamma=B$,
since by Lemma~\ref{lem:equivTT} (and its proof), if $P\in
TT(m/n,k,A)$ has data $d(P)=((N_r),\pi,(A,B,C))$, then $k>0$, $N_0=2$,
and $\pi(0,0)=(k,1)$: thus $TT(m/n,k,A)$ is empty if
$k\not\in\cA_{m/n}$, and has at most one element if
$k\in\cA_{m/n}$. The details are left to the reader.

\begin{thm}
\label{thm:caseA}
Let $0<i<\alpha$. Then $TT(m/n,k_i,A)$ consists of a single periodic
orbit~$P$, whose data $d(P)=((N_r),\pi,(A,B,C))$ satisfies: $N_r=2$
for $r\in T$, $N_r=1$ for $r\not\in T$, $A=\{0\}$, $B=\emptyset$,
$C=\{1\}$, and
\begin{eqnarray*}
\pi(0,0)&=&(k_i,1)\\
\pi(0,1)&=&(m,0)\\
\pi(r,s)&=&(r+_nm,s)\text{ if }r>0.
\end{eqnarray*}
$P$ has height $u_i/v_i$ and decoration $w_{m/n}$.
\end{thm}

\begin{exmp}
\label{ex:ttaeg}
The unique element of $TT(2/5,1,A)$ is the $\ast$-orbit of
Example~\ref{ex:dataeg}. 
\end{exmp}

\noindent
Combining the results of
Theorems~\ref{thm:hstt},~\ref{thm:casem1},~\ref{thm:identifym1},~\ref{thm:findperm},~\ref{thm:identify},~\ref{thm:caseA},~\ref{thm:totord}
and~\ref{thm:hs}b) gives the main result of this paper:
\begin{thm}
\label{thm:main}
Let $m/n\in(0,1/2]$. Then $TT(m/n)$ consists exactly of one orbit
$P_q^{m/n}$ of height $q$ and decoration $w_{m/n}$ for each rational
$q\in(0,m/n)$. The set of braid types of elements of $TT(m/n)$ is
totally ordered by the forcing relation, with
\mbox{$\bt(P_q^{m/n})\le\bt(P_{q'}^{m/n})$} if and only if $q\ge q'$.
\end{thm}

\section{Stars and the full horseshoe}
\label{sec:prune}
The aim of this section is to clarify the relationship between the
thick tree maps $F_{m/n}\colon T_n\to T_n$ and the full horseshoe
$F_{1/2}\colon T_2\to T_2$. The strongest connection is provided by
{\em pruning theory}. Each thick tree map $F_{m/n}$ can be obtained
from the full horseshoe by pruning: that is, by performing an isotopy
which destroys all of the dynamics of the horseshoe in an open subset
$U$ of $D^2$ (after the isotopy every point of $U$ is wandering),
while leaving the dynamics unchanged elsewhere (the isotopy is
supported on $U$). Although this process is conceptually valuable, it
is not required in the main body of the paper: since it is rather
complicated and assumes an understanding of the methods and results
of~\cite{dC,dCH1}, the treatment given in Section~\ref{sec:pfh} below
is on an intuitive level.

The only aspect of the relationship between stars and the full
horseshoe which was used in Section~\ref{sec:perorb} is the fact that
every $\ast$-orbit of $f_{m/n}$ has the same braid type as some
horseshoe periodic orbit. A rather straightforward proof of this is
given in Section~\ref{sec:hscode}: it has the additional advantage of
providing, for each $\ast$-orbit $P$, the code of a horseshoe periodic
orbit of the same braid type as $P$.

\subsection{The pruning approach}
\label{sec:pfh}
As stated above, the approach taken throughout this subsection is
intuitive: the interested reader is referred to~\cite{dCH1} for
details of the constructions outlined.

The aim is to start with the full horseshoe map $F_{1/2}$, to
perform a sequence of isotopies which decrease the dynamics
monotonically, remaining within the category of thick tree maps, and
to arrive at a given thick tree map $F_{m/n}\colon T_n\to T_n$. Two
distinct types of operation are used. The first redefines the thick
tree structure (that is, the subset of $D^2$ which is regarded as the
thick tree and the decomposition elements which give it its
structure), but leaves the dynamics unchanged. The second is an
isotopy which destroys the dynamics in some region which corresponds
to an interval in the underlying tree endomorphism whose image
backtracks. 

Both of these operations have counterparts on the level of the
underlying tree maps, and will be described on this level for the sake
of both conceptual and diagrammatic simplicity. The important point is
that the operations performed on the tree endomorphisms can be
realised by isotopies of the corresponding homeomorphisms of the
disk. Both the trees and their endomorphisms have additional structure
due to the fact that they are induced by thick tree maps: in
particular, the edges incident on a vertex have a cyclic order; and
when several edge images backtrack over a common edge, there is a
well-defined notion of the `innermost' backtracking.  This observation
is, of course, reflected in the way that tree maps have been drawn
throughout this paper.

Suppose that the image $f(e_r)=\ldots\be_se_s\ldots$ of an edge $e_r$
contains an innermost backtracking over the edge $e_s$. The following
two operations can be performed:
\begin{description}
\item[Glueing $\mathbf{e_r}$] Let $I$ be the subinterval of $e_r$ whose image is
$\be_se_s$. Identify those points of $I$ which have the same image (so
all points of $I$ except the preimage of the initial point of $e_s$
are identified in pairs). In the constructions described below, $I$
always has the central vertex $v$ of the star as an endpoint: glueing
therefore preserves the star structure, but increases the valence of
$v$ by $1$. It does not change the dynamics of $f$.
\item[Pulling tight $\mathbf{e_r}$] Delete the word $\be_se_s$ from $f(e_r)$. This
decreases the dynamics of $f$, while leaving the tree unchanged.
\end{description}

A second family of endomorphisms $g_{m/n}$ of $\Gamma_n$ are needed
for the construction: they differ from $f_{m/n}$ only in the image of
the edge $e_{n-1}$.

\begin{figure}[htbp]
\begin{center}
\pichere{0.6}{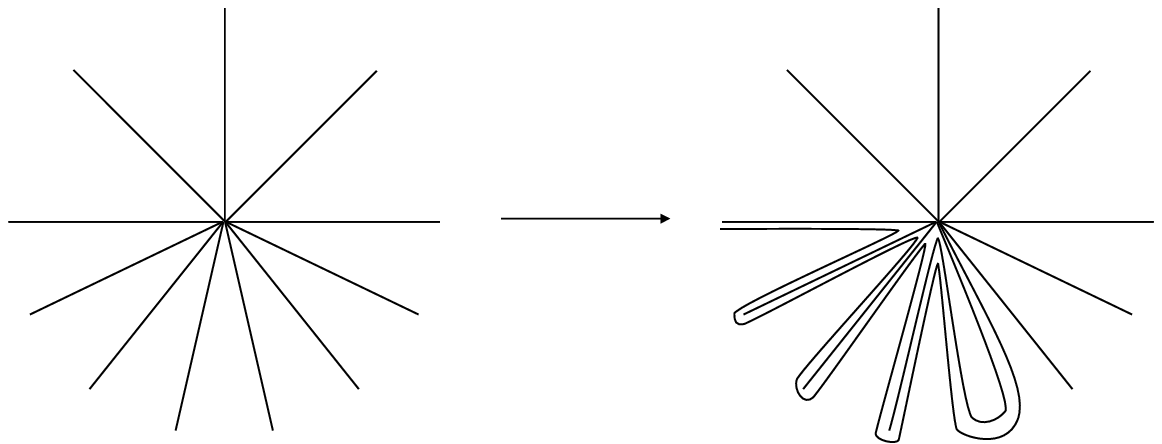}
\end{center}
\caption{The tree map $g_{4/11}\colon\Gamma_{11}\to\Gamma_{11}$}
\label{fig:gmn}
\end{figure}

Given $m/n\in(0,1/2)$, let $g_{m/n}\colon\Gamma_n\to\Gamma_n$ be the
tree map defined by
\begin{eqnarray*}
g(e_0)&=&e_0\be_1e_1\be_2e_2\ldots\be_me_m\\
g(e_r)&=&e_{r+_nm}\qquad(1\le r\le n-2)\\
g(e_{n-1})&=&e_{m-1}\be_me_m.
\end{eqnarray*}
As an example, $g_{4/11}$ is depicted in Fig.~\ref{fig:gmn}. These
maps are endowed with additional `2-dimensional' structure as
indicated in the figure.

The aim, then, is to construct the entire family $f_{m/n}$ of star
maps starting from $f_{1/2}$ by applying the operations of glueing and
pulling tight. This can be accomplished by a combination of two
compound procedures (together with the observation that one can pass
from $g_{m/n}$ to $f_{m/n}$ by pulling tight $e_{n-1}$). 
\begin{description}
\item[Procedure L] starts with $f_{m/n}$ and yields $g_{p/q}$, where
$p/q$ is the immediate left Farey child of $m/n$.
\item[Procedure R] starts with $g_{m/n}$ and yields $g_{u/v}$, where
$u/v$ is the immediate right Farey child of $m/n$.
\end{description}

To construct $f_{m/n}$ one navigates through the Farey graph from
$1/2$ to $m/n$, passing at each step from parent to immediate child:
at each step to the left (respectively right) one applies Procedure L
(respectively Procedure R). Thus, for example, to obtain $f_{3/10}$
from $f_{1/2}$ one applies Procedure L to obtain $g_{1/3}$; pulls
tight to obtain $f_{1/3}$; applies Procedure L to obtain $g_{1/4}$;
applies Procedure R to obtain $g_{2/7}$; applies Procedure R to obtain
$g_{3/10}$; and pulls tight to obtain $f_{3/10}$.

\par\noindent{\bf Procedure L: }\ 
Start with $f_{m/n}$ and glue $e_0$ (that is, identify pairs of points
in $e_0$ which have the same image in $e_m$). This creates a new edge,
which is labelled $e_n$. Now the image of $e_{n-m}$ backtracks over
$e_n$; glue $e_{n-m}$. This creates another new edge $e_{n+1}$, and
the image of $e_{n-2m}$ backtracks over it. Continue this procedure
until edge $e_{m-1}$ has been glued; each time an edge $e_j$ is glued
with $1\le j\le m-1$, pull tight the innermost backtracking of
the image of $e_0$ over $e_j$ before proceeding. This yields
$g_{p/q}$, where $p/q$ is the immediate left Farey child of $m/n$.

Applying Lemma~\ref{lem:fareyprops}, observe that if $\LFP(m/n)=a/b$,
a total of $\card\cO_{m/n}[m-1,0]=b$ glueings are performed, and hence
at the end of the procedure the star has $n+b=q$ edges; and
$\card(\cO_{m/n}[m-1,0]\cap[0,m-1])=a$ of the edges glued are between
$e_0$ and $e_m$ so that, after relabelling the edges in cyclic order,
the image of $e_0$ crosses edges $e_0$ to $e_{m+a}=e_p$.

\begin{exmp}
Let $m/n=3/7$: thus procedure~L gives a construction of $g_{5/12}$ from
 $f_{3/7}$. The edges $e_0$, $e_4$, $e_1$, $e_5$, and $e_2$
are glued successively: the new edges thereby created are labelled
$e_7$, $e_8$, $e_9$, $e_{10}$ and $e_{11}$ respectively. Since
$m-1=2$, the only pulling tights occur after glueing $e_1$ and $e_2$.
The procedure is shown in Fig.~\ref{fig:glueitup}. 

\begin{figure}[htbp]
\Lab{0}{0}{r}
\Lab{1}{1}{tr}
\Lab{2}{2}{t}
\Lab{3}{3}{tl}
\Lab{4}{4}{l}
\Lab{5}{5}{bl}
\Lab{6}{6}{br}
\Lab{7}{7}{r}
\Lab{8}{8}{l}
\Lab{9}{9}{tr}
\Lab{a}{10}{b}
\Lab{b}{11}{tl}
\Lab{d}{f(5)=1}{tr}
\Lab{e}{f(6)=2}{t}
\Lab{f}{f(1)=4}{l}
\Lab{g}{f(2)=5}{bl}
\Lab{h}{f(3)=6}{}
\Lab{i}{f(4)=0}{}
\psfrag{j}[b]{{\tiny\bf Glue} $\scriptscriptstyle{\mathbf{0}}$}
\Lab{k}{f(7)=3}{tl}
\psfrag{l}[l]{{\tiny\bf Glue} $\scriptscriptstyle{\mathbf{4}}$}
\Lab{m}{f(8)=7}{tr}
\psfrag{n}[b]{{\tiny\bf Glue} $\scriptscriptstyle{\mathbf{1}}$}
\Lab{o}{f(9)=8}{bl}
\psfrag{p}[l]{{\tiny\bf Pull Tight}}
\psfrag{q}[b]{{\tiny\bf Glue} $\scriptscriptstyle{\mathbf{5}}$}
\Lab{r}{f(10)=9}{tr}
\psfrag{s}[l]{{\tiny\bf Glue} $\scriptscriptstyle{\mathbf{2}}$}
\Lab{t}{f(11)=10}{b}
\psfrag{u}[b]{{\tiny\bf Pull Tight}}
\begin{center}
\pichere{0.9}{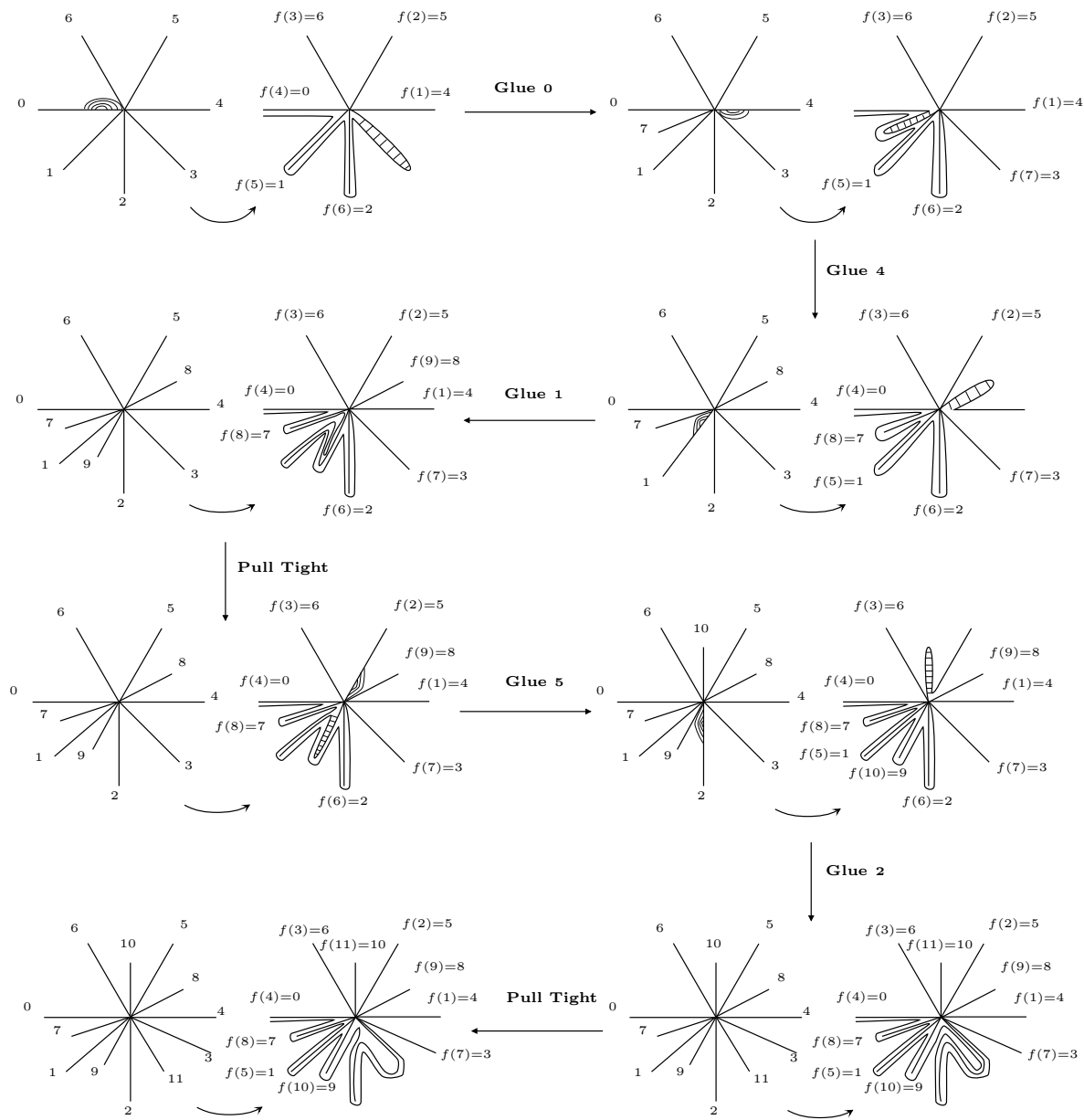}
\end{center}
\caption{Constructing $g_{5/12}$ from $f_{3/7}$}
\label{fig:glueitup}
\end{figure}

\end{exmp}

\par\noindent{\bf Procedure R: }\ is similar. Start with
 $g_{m/n}$. Since $g_{m/n}(e_{n-1})$ backtracks over $e_m$ one can
 successively glue $e_{n-1}$, $e_{n-m-1}$, $e_{n-2m-1}$, \ldots,
 $e_m$, pulling tight $e_0$ each time an edge $e_j$ is glued with
 $1\le j\le m-1$. This yields $g_{u/v}$, where $u/v$ is the immediate
 right Farey child of $m/n$.

\subsection{Horseshoe symbolics for $\ast$-orbits}
\label{sec:hscode}
In this section a more direct approach to the problem of showing that
$\ast$-orbits have horseshoe braid type is outlined. It makes use of
the notion of the {\em line diagram} of an isotopy class of
homeomorphisms of the punctured disk. For each $n\ge 2$, let $D_n$ be
a standard model of the $n$-punctured disk in which the punctures (or
marked points) are equally spaced along the horizontal diameter of the
disk. The line diagram of a homeomorphism $f\colon D_n\to D_n$ is the
sequence $([f(\alpha_1)], \ldots,[f(\alpha_{n-1})])$ of homotopy
classes of the images of the horizontal arcs $\alpha_i$ joining the
$i^{\mbox{\scriptsize{th}}}$ to the $(i+1)^{\mbox{\scriptsize{th}}}$
puncture. It is straightforward to show that two homeomorphisms
$f,g\colon D_n\to D_n$ are isotopic if and only if they have the same
line diagram.

If $F\colon D^2\to D^2$ is an orientation-preserving homeomorphism
having a period $n$ orbit~$P$, then conjugating $F$ so that the points
of $P$ coincide with the punctures of $D_n$ yields a homeomorphism of
$D_n$ and thence a line diagram. Different conjugacies naturally give
rise to different line diagrams, but periodic orbits $P$ and $Q$ of
homeomorphisms $F$ and $G$ have the same braid type if and only if the
conjugacies can be chosen so as to give the same line diagrams.

If $F$ is the horseshoe, then there is a natural choice of (isotopy
class of) conjugacy, namely one which sends vertical leaves to
vertical leaves in an order-preserving manner. Such a conjugacy gives
rise to `unimodal' line diagrams which reflect the underlying unimodal
structure of the horseshoe: for example, the periodic orbit of code
$10010$ gives rise to the line diagram of Fig.~\ref{fig:simporb}.

Thus to show that a $\ast$-orbit $P$ of $F_{m/n}$ has horseshoe braid
type, it is necessary to construct a conjugacy which places the points
of $P$ along the horizontal diameter of $D^2$ in such a way that the
resulting line diagram is unimodal. This will be achieved by drawing
an arc through the points of the orbit which will be mapped by the
conjugacy onto an interval of the horizontal diameter of $D^2$. Not
only does this show that $P$ has horseshoe braid type, it also yields
the code of a horseshoe periodic orbit of the same braid type as $P$:
those points whose images are on the increasing segment of the
unimodal line diagram are coded $0$, while those on the
decreasing segment are coded $1$ -- the critical point can be coded
either $0$ or $1$.

\begin{figure}[htbp]
\lab{a}{f_{2/5}}{b}
\begin{center}
\pichere{0.6}{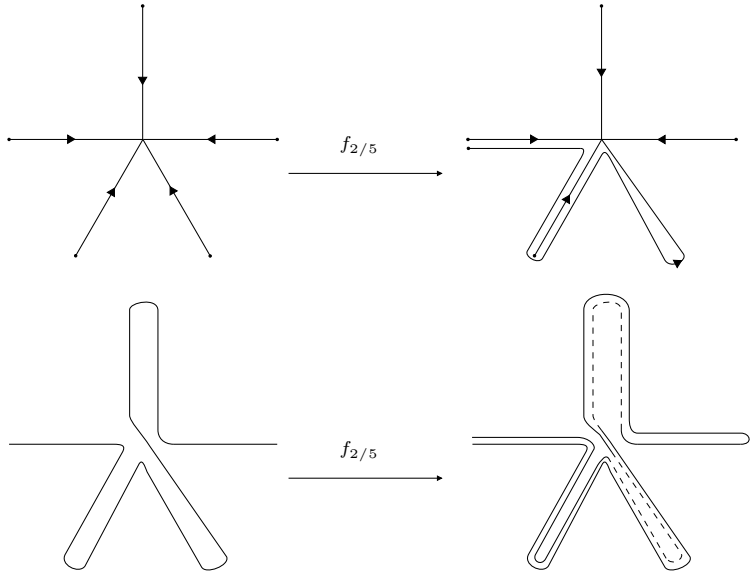}
\end{center}
\caption{A path whose image under $f_{2/5}$ is a subhorseshoe }
\label{fig:hsarc}
\end{figure}

The crucial observation is that if $\alpha_{m/n}$ denotes the path
\[\alpha_{m/n}=e_0\be_1e_1\be_2e_2\ldots\be_{n-m-1}e_{n-m-1}\be_{n-1}e_{n-1}\be_{n-2}e_{n-2}\ldots
\be_{n-m+1}e_{n-m+1}\be_{n-m},\] then its image 
\[f_{m/n}(\alpha_{m/n})=e_0\be_1e_1\ldots\be_{n-1}e_{n-1}\be_{m-1}e_{m-1}\be_{m-2}e_{m-2}\ldots\be_1e_1\be_0\]
can be obtained from the image
\[h_{m/n}(\alpha_{m/n})=\alpha_{m/n}\overline{\alpha}_{m/n}=
e_0\be_1e_1\be_2e_2\ldots\be_{n-m-1}e_{n-m-1}\]\[\left(
\be_{n-1}e_{n-1}\be_{n-2}e_{n-2}\ldots\be_{n-m+1}e_{n-m+1}\right)
\be_{n-m}e_{n-m}\be_{n-m+1}\ldots\be_{n-1}e_{n-1}\]\[\left(
\be_{n-m-1}e_{n-m-1}\ldots\be_me_m\right)
\be_{m-1}e_{m-1}\ldots\be_1e_1\be_0\] 
of $\alpha_{m/n}$ under the horseshoe $h_{m/n}$ by removing the
bracketed words (see Fig.~\ref{fig:hsarc}, which illustrates this
for $m/n=2/5$). Notice that if $\Gamma_n$ is drawn (as in this paper)
so that $e_0$ and $e_{n-m}$ are horizontal, then $\alpha_{m/n}$ passes
along $e_0$, then around each of the edges below the horizontal in the
positive direction, then around each of the edges above the horizontal
in the negative direction, and finally along $e_{n-m}$.

It follows that if $P$ is a $\ast$-orbit of $f_{m/n}$, then $P$ has
horseshoe braid type provided that an arc $\alpha$ projecting to
$\alpha_{m/n}$ can be passed through the points of $P$ in the thick
tree $T_n$ in such a way that $F_{m/n}(\alpha)$ projects to
$f_{m/n}(\alpha_{m/n})$. Moreover, horseshoe symbolics for $P$ can be
obtained by coding with $0$ those points of $P$ which either lie in
$e_0\cup\cdots\cup e_{n-2m-1}$, or lie in $e_{n-2m}$ and have $\alpha$
pass through them with the orientation of $\be_{n-2m}$ (these are
precisely the points whose images lie in the increasing segment of the
unimodal line diagram); and coding with $1$ the other points of $P$.

The only issue in constructing such an arc is to decide whether it
should pass through each point of $P$ with the orientation of $e_r$,
or with the orientation of $\be_r$. The following partition of
$\cL\setminus\{(r,0):0\le r<n\}$ gives the unique coherent way of
doing this: points which $\alpha$ passes through with the orientation
of $e_r$ lie in $I$.

\begin{defn}
Let $P\in\cP(m/n)$ have data $d(P)=((N_r),\pi,(A,B,C))$. Define a
partition 
\[\cL\setminus\{(r,0):0\le r<n\}=I\cup O\] inductively as
follows. $\pi(0,0)$ lies in $I$ (respectively $O$) if $0\in B$
(respectively $0\in A$). For each $i$ with $2\le i\le\card P-n$, let
$\pi^{i-1}(0,0)=(r,s)$. Then $\pi^i(0,0)$ and $\pi^{i-1}(0,0)$ are in
the same set if $r=0$ and $s\in B$, or if $0<r<n-2m$. They are in
different sets if $r=0$ and $s\in A$, or if $n-2m<r<n$. If $r=n-2m$,
then $\pi^i(0,0)\in O$.
\end{defn}

The following theorem then follows from the discussion above.

\begin{thm}
\label{thm:hscode}
Let $P$ be a period $N$ orbit in $\cP(m/n)$ for some $m/n\not=1/2$.
Define a map $c\colon\cL\to\{0,1\}$ by
\[c(r,s)=\left\{\begin{array}{ll}
0 & \mbox{ if } r<n-2m \mbox{ or }r=n-2m\mbox{ and }(r,s)\in O\mbox{
or }(r,s)=(n-2m,0)\\
1 & \mbox{ if }r>n-2m\mbox{ or }r=n-2m\mbox{ and }(r,s)\in I.
\end{array}\right.\]

Then $P$ has the same braid type as the horseshoe periodic orbit of
code
\[c(n-m,0)c(\pi(n-m,0))c(\pi^2(n-m,0))\ldots c(\pi^{N-1}(n-m,0)).\]
\end{thm}

\begin{exmp}
Let $P\in\cP(2/5,1,A)$ be the periodic orbit with data
\[
\ovun{A}{I}{(0,0)}\to
\un{O}{(1,1)}\to
\un{O}{(3,1)}\to
\ovun{C}{I}{(0,1)}
\to(2,0)\to(4,0)\to(1,0)\to(3,0).
\]
The partition $\cL\setminus\{(r,0)\colon r>0\}=I\cup O$ is indicated
with this data. The horseshoe symbolics of $P$ can then be
written down:
\[
\un{0}{(0,0)}\to
\un{0}{(1,1)}\to
\un{1}{(3,1)}\to
\un{0}{(0,1)}\to
\un{1}{(2,0)}\to
\un{1}{(4,0)}\to
\un{0}{(1,0)}\to
\un{1}{(3,0)}.
\]

Hence $P$ has the same braid type as the horseshoe orbit of code
$10010110$. Thus $P$ has height $1/3$ and decoration $11=w_{2/5}$, as
expected since $P$ is the unique element of $TT(2/5,1,A)$ discussed in
examples~\ref{ex:dataeg} and~\ref{ex:ttaeg}.
\end{exmp}

Notice that since the cycle notation of $\pi$ always ends
$(m,0)\to(2m,0)\to\cdots\to(n-m,0)$, every element of $\cP(m/n)$ has
the same braid type as a horseshoe orbit whose code ends $w_{m/n}0$.

\section{Proof of Lemma~\ref{lem:identify}}
\label{sec:app}

The following lemma summarizes the results of~\cite{Ha} which will be
used in the proof of Lemma~\ref{lem:identify}. Part~a) is theorem~3.10
of~\cite{Ha}, part~b) is theorem~3.11, part~c) is lemma~3.4 (in which
the notation $d_{r/s}$ is used to mean $0w_{r/s}$), and part~d) is a
combination of theorem~3.5 and lemma~3.6. The first part of the lemma
gives an algorithm for computing the rotation interval of an arbitrary
horeseshoe periodic orbit. Although it is complicated to state, it is
much easier to explain intuitively. For each block of $0$s in $c_P$,
calculate the heights of the sequences obtained by moving forwards and
backwards through $\overline{c_P}$ starting at the $1$ immediately
before the block of $0$s: if the backward sequence starts $11$, then
first replace the second $1$ with a $0$. If the backward height is not
less than the forward height, then the interval between the two is
contained in the rotation interval. The union of all such intervals is
the rotation interval.

\begin{lem}
\label{lem:riprops}
\begin{enumerate}[a)]
\item Let $P$ be a horseshoe
periodic orbit which contains the point of itinerary
\[\overline{0^{\kappa_1}1^{\mu_1}0^{\kappa_2}1^{\mu_2}\ldots
0^{\kappa_r} 1^{\mu_r}}\]
(written in such a way that $\kappa_i,\mu_i>0$ for all $i$). Then $P$
has rotation interval
\[\rho i(P)=\bigcup_{i=1}^r[\xi_i,\eta_i],\]
where
\[\xi_i=q\left(1\left(
\overline{0^{\kappa_i}1^{\mu_i}0^{\kappa_{i+1}}1^{\mu_{i+1}}\ldots
1^{\mu_r}0^{\kappa_1}\ldots0^{\kappa_{i-1}}1^{\mu_{i-1}}}
\right)\right)\qquad\text{ and}\]
\[\eta_i=\left\{\begin{array}{ll}
q\left(\overline{10^{\kappa_{i-1}}1^{\mu_{i-2}}0^{\kappa_{i-2}}\ldots
0^{\kappa_1}1^{\mu_r}\ldots 1^{\mu_i}0^{\kappa_i}}\right)
&\quad\text{ if }\mu_{i-1}=1\\
q\left(10\left(\overline{
1^{\mu_{i-1}-2}0^{\kappa_{i-1}}1^{\mu_{i-2}}0^{\kappa_{i-2}}\ldots
0^{\kappa_1}1^{\mu_r}\ldots 1^{\mu_i}0^{\kappa_i}1^2
}\right)\right)
&\quad\text{ if }\mu_{i-1}>1
\end{array}\right.\]
(and $[\xi_i,\eta_i]=\emptyset$ if $\eta_i<\xi_i$).
\item Let $P$ be a horseshoe periodic orbit. Then the left hand
endpoint of $\rho i(P)$ is $q(P)$.
\item Let $0<r/s<1/2$. Then $c\in\{0,1\}^\Nset$ has height $q(c)=r/s$ 
if and only if
\[\overline{10w_{r/s}1} \preceq c \preceq 10\,\,
\overline{w_{r/s}011}\]
(where the inequalities are with respect to the unimodal order on
$\{0,1\}^\Nset$).  In particular, if $q(c)=r/s$ then
$c=10w_{r/s}\opti1\ldots$, and the first isolated $1$ in $c$ cannot
appear before the $s+1^{\mbox{\scriptsize{th}}}$ symbol.
\item If $P$ is a period $N$ horseshoe orbit which is not of finite
order braid type and which has height $q(P)=u/v$, then $N\ge
v+2$. Moreover, if $N=v+2$ then $c_P=c_{u/v}\opti$ (and $P$ has
rotation interval $[u/v,1/2]$).
\end{enumerate}
\end{lem}

The following lemma will also be needed:
\begin{lem}
\label{lem:tweak}
Let $m/n<1/2$, and write
$c_{m/n}=10^{\kappa_1}1^20^{\kappa_2}1^2\ldots 1^20^{\kappa_m}1$
(where the $\kappa_i$ are given by formula~(\ref{eq:kappa}) on
page~\pageref{eq:kappa}). Let $1\le r\le m$. Then the word
\[c=10^{\kappa_r+1}1^20^{\kappa_{r+1}}1^2\ldots 1^20^{\kappa_m}1\]
disagrees with $c_{m/n}$ within the shorter of their lengths,
and is greater than it in the unimodal order.
\end{lem}

\begin{proof}
If the two words didn't disagree, then it would follow that
$\kappa_r+1=\kappa_1$ and that $\kappa_m=\kappa_{m-r+1}$,
contradicting the fact that $c_{m/n}$ is palindromic.

Observe that formula~(\ref{eq:kappa}) gives, for each $s$ with
$1\le s\le m+1-r$,
\[\sum_{i=1}^s\kappa_i=\left\lfloor\frac{sn}{m}\right\rfloor-(2s-1),\]
and
\begin{eqnarray*}
\kappa_r+1+\sum_{i=r+1}^{r+s-1}\kappa_i
&=&\left\lfloor\frac{(r+s-1)n}{m}\right\rfloor-\left\lfloor
   \frac{(r-1)n}{m}\right\rfloor-2s+1\\
&\ge&\left\lfloor \frac{(r+s-1)n}{m}-\frac{(r-1)n}{m}\right\rfloor
   -2s+1\\
&=&\left\lfloor\frac{sn}{m}\right\rfloor-(2s-1).
\end{eqnarray*}
Hence at the point where they first disagree $c$ has a longer block
of $0$s than $c_{m/n}$, and so is greater in the unimodal order.
\end{proof}

\noindent{\bf Lemma~\ref{lem:identify}}\,\,{\em Let $P$ be a period
$N>1$ orbit of the horseshoe with non-trivial rotation interval
$\rho i(P)=[u/v,m/n]$. Then $N\ge v+n$. Moreover, $N=v+n$ if and only
if $c_P$ is one of the four words $c_{u/v}\opti w_{m/n}\opti$ (or one
of the two words $c_{u/v}\opti$ in the case $m/n=1/2$).}

\medskip

\begin{proof}
By induction on $N$, with the case $N=2$ vacuous since the only period
$2$ horseshoe orbit has trivial rotation interval. The case $m/n=1/2$
follows immediately from parts~b) and~d) of Lemma~\ref{lem:riprops},
so it will be assumed that $m/n<1/2$. By parts~b) and~d) of
Lemma~\ref{lem:riprops}, it then follows that $N\ge v+3$: let
$k=N-v\ge 3$. By Theorem~\ref{thm:hs}~a), the code of $P$ is of the
form
\[c_P=c_{u/v}\opti w\opti,\]
for some word $w$ of length $k-3$. The four cases $c_P=c_{u/v}0w0$,
$c_{u/v}0w1$, $c_{u/v}1w0$, and $c_{u/v}1w1$ need to be considered
separately. 

\noindent{\bf Case a)} Suppose that
\begin{eqnarray*}
c_P&=&c_{u/v}0w0\\
&=& 10^{\kappa_1}1^20^{\kappa_2}1^2\ldots 1^20^{\kappa_u}\,\,1
0^{\lambda_1}1^{\nu_1}0^{\lambda_2}1^{\nu_2}\ldots 1^{\nu_{t-1}}0^{\lambda_t},
\end{eqnarray*}
where $\kappa_i=\kappa_i(u/v)$ for $1\le i\le u$, and $\lambda_i$
and $\nu_i$ are chosen to be positive for all $i$. Note that the
initial word of $c_P$ up to $0^{\kappa_u}$ has length $v$, and the
complementary word has length $k$.
For $1\le i\le u$,
let $[\xi_i^1,\eta_i^1]$ be the contribution to the rotation interval
corresponding to the block $0^{\kappa_i}$ of $0$s; and for $1\le i\le
t$, let $[\xi_i^2,\eta_i^2]$ be the contribution corresponding to
$0^{\lambda_i}$. Thus
\[\rho i(P)=\bigcup_{i=1}^u\left[\xi_i^1,\eta_i^1\right]\,\,\cup
\,\,\bigcup_{i=1}^t\left[\xi_i^2,\eta_i^2\right].\]

Now for $2\le i\le u$,
$\eta_i^1=q(10^{\kappa_{i-1}+1}1^20^{\kappa_{i-2}}1^2\ldots
1^20^{\kappa_1}10\ldots)$, which, by Lemmas~\ref{lem:tweak}
and~\ref{lem:riprops}~c) and the fact that $c_{u/v}$ is palindromic,
is less than $u/v$. Hence, since $\rho i(P)=[u/v,m/n]$, the interval
$[\xi_i^1,\eta_i^1]$ must be empty.  On the other hand, $\xi_1^1=u/v$,
and
\[\eta_1^1=q\left(10^{\lambda_t}1^{\mu_{t-1}}\ldots
1^{\mu_1}0^{\lambda_1}10\ldots\right)\]
Let $\eta_1^1=r/s$: then by Lemma~\ref{lem:riprops}~c), $k+1\ge s+1$,
or $k\ge s$. It follows that $\bigcup_{i=1}^u\left[\xi_i^1,\eta_i^1\right]$
is equal either to $\{u/v\}$, or to $[u/v,r/s]$ for some $r/s$ with
$s\le k=N-v$. 

Now consider the intervals $[\xi_i^2,\eta_i^2]$. First,
$\eta_1^2=u/v$, while $\xi_1^2\ge u/v$ (since $\overline{c_P}$ has
height $u/v$ and $q\colon\{0,1\}^\Nset\to(0,1/2]$ is
order-reversing). Hence $[\xi_1^2,\eta_1^2]$ is either empty or equal
to $\{u/v\}$. By Lemma~\ref{lem:riprops}~c), all of the other
$[\xi_j^2,\eta_j^2]$ are either empty or of the form $[a/b,c/d]$,
where both $b$ and $d$ are less than $k$. 

Hence $\rho i(P)=[u/v,m/n]$, where $n\le k=N-v$, and the first statement
of the lemma follows. For the second statement, observe from the
argument above that if $n=k$ then 
\[m/n=\eta_1^1=q\left(10^{\lambda_t}1^{\mu_{t-1}}\ldots
1^{\mu_1}0^{\lambda_1}10\ldots\right),\]
and using Lemma~\ref{lem:riprops}~c) it follows that $w=w_{m/n}$ as
required.

\noindent{\bf Case b)}
Suppose that $c_P=c_{u/v}0w1$ for some word $w$ of length
$k-3$. If $c_{u/v}0w0$ is a maximal (non-repetitive) word, then it is
the code of an orbit of the same braid type as $P$, and the result
follows from the case~a). If it is not maximal, then $P$ is the
period-doubling of a horseshoe orbit of half its period with the same
rotation interval, and the result follows from the inductive
hypothesis.

\noindent{\bf Case c)} If
$c_P=c_{u/v}1w0=10^{\kappa_1}1^20^{\kappa_2}1^2\ldots
1^20^{\kappa_u}\,\,1^{\nu_1}
0^{\lambda_1}1^{\nu_2}0^{\lambda_2}1^{\nu_3}\ldots
1^{\nu_t}0^{\lambda_t}$, decompose $\rho i(P)$ into intervals
$[\xi_i^1,\eta_i^1]$ for $1\le i\le u$ and $[\xi_i^2,\eta_i^2]$ for
$1\le i\le t$ as in  case~a). Just as in that case, it can be
shown that $[\xi_i^1,\eta_i^1]$ is empty for $2\le i\le u$. Thus
either $\eta_1^1$ or some $\eta_i^2$ is equal to $m/n$. By
Lemma~\ref{lem:riprops}~a) and~c), this means that the reverse
$\widehat{c_P}$ of the code of $P$ has the property that
$\overline{\widehat{c_P}}$ contains one of the words $01\opti
w_{m/n}\opti 1$: equivalently (since $w_{m/n}$ is palindromic),
$\overline{c_P}$ contains one of the words $1\opti w_{m/n}\opti10$. If
there is such a word which is disjoint from the prefix
$10^{\kappa_1}1^20^{\kappa_2}1^2\ldots 1^20^{\kappa_u}1$ of $c_P$,
then the proof can be completed as in the case~a). It remains to
show, therefore, that if such a word overlaps the prefix, then the
corresponding contribution to $\rho i(P)$ can not have right hand
endpoint $m/n$. (It is not necessary to consider words contained
entirely within the prefix, since it has already been shown that
$\eta_i^1<u/v$ for $i\ge 2$.)

Write $\kappa'_i=\kappa_i(m/n)$ for $1\le i\le m$, and suppose first
that $\overline{\widehat{c_P}}$ contains a word
\[010w_{m/n}01=010^{\kappa_1'}1^20^{\kappa_2'}1^2\ldots1^2
0^{\kappa_m'}1\] whose associated interval has right hand endpoint
$m/n$, and suppose that the final block $0^{\kappa_m'}$ of $0$s in
this (reverse) word coincides with the block $0^{\kappa_r}$ of $0$s in
the prefix. Then $r\ge2$, since otherwise $P$ would have height
$m/n$. By Lemma~\ref{lem:riprops}~a), the right hand
endpoint of the associated interval being $m/n$ gives
\[q(10^{\kappa_1'}1^2\ldots1^20^{\kappa_m'}1^20^{\kappa_{r-1}}1^2\ldots
1^20^{\kappa_1}10\ldots)=q(10w_{m/n}0110^{\kappa_{r-1}}1^2\ldots1^20^{\kappa_1}10\ldots)=m/n.\]
Hence by Lemma~\ref{lem:riprops}~c), $0^{\kappa_{r-1}}1^2
\ldots1^20^{\kappa_1}10\ldots\succ\overline{w_{m/n}011}$, or
equivalently
\[10^{\kappa_{r-1}+1}1^20^{\kappa_{r-2}}1^2\ldots1^20^{\kappa_1}10\ldots
\prec 10\overline{w_{m/n}011}.\]
Applying Lemma~\ref{lem:riprops}~c) again gives
$q(10^{\kappa_{r-1}+1}1^20^{\kappa_{r-2}}1^2\ldots1^20^{\kappa_1}10\ldots)
\ge m/n>u/v$, so that
\[10^{\kappa_{r-1}+1}1^20^{\kappa_{r-2}}1^2\ldots1^20^{\kappa_1}10\ldots
\prec 10^{\kappa_1}1^20^{\kappa_2}1^2\ldots0^{\kappa_m}10\ldots,\]
contradicting Lemma~\ref{lem:tweak}.

Exactly the same argument works if the contributing word in
$\overline{\widehat{c_P}}$ is $011w_{m/n}01$. If the word is one of
$01\opti w_{m/n}11$, then comparing the overlapping segments of this
word and the prefix of $P$ gives $\kappa_r=\kappa'_m-1$ and hence
(using again that $w_{m/n}$ is palindromic)
\begin{eqnarray*}
10^{\kappa_1'}1^20^{\kappa_2'}1^2\ldots 1^20^{\kappa_{m-r}'}
1^20^{\kappa_{m-r+1}'}\ldots 
&=&10^{\kappa_r+1}1^20^{\kappa_{r+1}}1^2\ldots 1^20^{\kappa_{m-1}}
1^20^{\kappa_m}\ldots\\
&\succ& 10^{\kappa_1}1^20^{\kappa_2}1^2\ldots1^20^{\kappa_m}\ldots,
\end{eqnarray*}
(where the final inequality is by Lemma~\ref{lem:tweak}). Taking
heights gives
\[m/n=q(10^{\kappa_1'}1^20^{\kappa_2'}1^2\ldots 1^20^{\kappa_{m-r}'}
1^20^{\kappa_{m-r+1}'}\ldots) <
q(10^{\kappa_1}1^20^{\kappa_2}1^2\ldots1^20^{\kappa_m}\ldots) =u/v,\]
a contradiction.

\noindent{\bf Case~d)} The proof for $c_P=c_{u/v}1w1$ follows from
case~c) in the same way that case~b) follows from case~a).
\end{proof}

\bibliography{star}

\providecommand{\bysame}{\leavevmode\hbox to3em{\hrulefill}\thinspace}
\begin{thebibliography}{BGMY80}

\bibitem[BGH93]{BGH}
P.~Boyland, J.~Guaschi, and T.~Hall, \emph{L'ensemble de rotation des
  hom\'eomorphismes pseudo-{A}nosov}, C. R. Acad. Sci. Paris S\'er. I Math.
  \textbf{316} (1993), no.~10, 1077--1080.

\bibitem[BGMY80]{BGMY}
L.~Block, J.~Guckenheimer, M.~Misiurewicz, and L.S. Young, \emph{Periodic
  points and topological entropy of one-dimensional maps}, Global theory of
  dynamical systems (Proc. Internat. Conf., Northwestern Univ., Evanston, Ill.,
  1979), Springer, Berlin, 1980, pp.~18--34.

\bibitem[BH95]{BH}
M.~Bestvina and M.~Handel, \emph{Train-tracks for surface homeomorphisms},
  Topology \textbf{34} (1995), no.~1, 109--140.

\bibitem[Boy84]{Boy}
P.~Boyland, \emph{Braid types and a topological method of proving positive
  entropy}, Preprint, Boston University, 1984.

\bibitem[Boy92]{Boyri}
P.~Boyland, \emph{Rotation sets and monotone periodic orbits for annulus
  homeomorphisms}, Comment. Math. Helv. \textbf{67} (1992), no.~2, 203--213.

\bibitem[Boy94]{Boy2}
P.~Boyland, \emph{Topological methods in surface dynamics}, Topology Appl.
  \textbf{58} (1994), no.~3, 223--298.

\bibitem[dC99]{dC}
A.~de~Carvalho, \emph{Pruning fronts and the formation of horseshoes}, Ergodic
  Theory Dynam. Systems \textbf{19} (1999), no.~4, 851--894.

\bibitem[dCHa]{invs}
A.~de~Carvalho and T.~Hall, \emph{Conjugacies between horseshoe braids}, In
  preparation.

\bibitem[dCHb]{dCH2}
A.~de~Carvalho and T.~Hall, \emph{The forcing relation for horseshoe braid
  types}, To appear in Experiment. Math.

\bibitem[dCH01]{dCH1}
A.~de~Carvalho and T.~Hall, \emph{Pruning theory and {T}hurston's
  classification of surface homeomorphisms}, J. Eur. Math. Soc. (JEMS)
  \textbf{3} (2001), no.~4, 287--333.

\bibitem[Hal]{programs}
T.~Hall, Software available from {\tt http://www.liv.ac.uk/\~{}tobyhall/hs/}.

\bibitem[Hal94]{Ha}
T.~Hall, \emph{The creation of horseshoes}, Nonlinearity \textbf{7} (1994),
  no.~3, 861--924.

\bibitem[Han90]{Handelrot}
M.~Handel, \emph{The rotation set of a homeomorphism of the annulus is closed},
  Comm. Math. Phys. \textbf{127} (1990), no.~2, 339--349.

\bibitem[HW79]{Farey}
G.~Hardy and E.~Wright, \emph{An introduction to the theory of numbers}, fifth
  ed., The Clarendon Press Oxford University Press, New York, 1979.

\bibitem[Sma67]{Smale}
S.~Smale, \emph{Differentiable dynamical systems}, Bull. Amer. Math. Soc.
  \textbf{73} (1967), 747--817.

\bibitem[Thu88]{Thu}
W.~Thurston, \emph{On the geometry and dynamics of diffeomorphisms of
  surfaces}, Bull. Amer. Math. Soc. (N.S.) \textbf{19} (1988), no.~2, 417--431.

\end{thebibliography}
\end{document}